# Analysis of the navigation of magnetic microrobots through cerebral bifurcations

Pedro G Alves,[1,2] Maria Pinto,[3] Rosa Moreira,[3] Derick Sivakumaran,[4] Fabian Landers,[5] Maria Guix,[7] Bradley J Nelson,[5] Andreas D Flouris,[6] Salvador Pané,[5] Josep Puigmartí-Luis,[7,8,*] Tiago S Mayor[1,2,*]

[1] Transport Phenomena Research Centre (CEFT), Engineering Faculty, Porto University, Portugal
[2] Associate Laboratory in Chemical Engineering (ALICE), Engineering Faculty, Porto University, Portugal
[3] Experian, Portugal
[4] Magnebotix, Zurich, Switzerland
[5] Multi-Scale Robotics Lab, ETH Zurich, Zurich, Switzerland
[6] FAME Laboratory, Department of Exercise Science, University of Thessaly, Greece
[7] Depart. de Ciència dels Materials i Química Física and Institut de Química Teòrica i Computacional, University of Barcelona, Barcelona, Spain
[8] Institució Catalana de Recerca i Estudis Avançats (ICREA), Barcelona, Spain

* Corresponding authors: josep.puigmarti@ub.edu; tiago.sottomayor@fe.up.pt

# Abstract

Local administration of thrombolytics in ischaemic stroke could accelerate clot lysis and the ensuing reperfusion while minimizing the side effects of systemic administration. Medical microrobots could be injected into the bloodstream and magnetically navigated to the clot for administering the drugs directly to the target. The magnetic manipulation required to navigate medical microrobots will depend on various parameters such as the microrobots size, the blood velocity, and the imposed magnetic field gradients. Numerical simulation was used to study the motion of magnetically controlled microrobots flowing through representative cerebral bifurcations, for predicting the magnetic gradients required to navigate the microrobots from the injection point until the target location. Upon thorough validation of the model against several independent analytical and experimental results, the model was used to generate maps and a predictive equation providing quantitative information on the required magnetic gradients, for different scenarios. The developed maps and predictive equation are crucial to inform the design, operation and optimization of magnetic navigation systems for healthcare applications.

## 1. Introduction

Stroke is characterized by an inadequate blood supply to cerebral tissue, which can lead to long-term disability and death. According to the 2019 Global Burden of Disease study, stroke is the second-leading cause of death worldwide and its incidence has increased by 70% between 1990 and 2019. Of all the worldwide new stroke occurrences in 2019, two-thirds (7.6 million cases) were of the ischaemic type.[1] Ischaemic strokes are characterized by a vascular occlusion that blocks blood flow to downstream arteries, leading to tissue necrosis and serious long-term complications including death.[2] The established treatment for ischaemic stroke is the systemic administration of thrombolytic agents to dissolve the clot obstructing the artery.[3,4] However, only up to 20% of cases respond to treatment as the clots can be too big to dissolve in a timely manner, due to the systemic administration instead of a localized approach.[5] Dosage increase is not possible due to the toxicity and side effects of this type of drugs, including blood thinning effects that can lead to bleeding.[2,6] Mechanical removal of the clots (or thrombectomy) is an additional approach where a catheter is used to reach into the cerebral arteries to extract the clot and improve recanalization.[5,7] Yet, that can be challenging in narrow and tortuous arteries, making it difficult to navigate the catheter to deep locations. Hence, a new approach to reach deep neurovascular regions and locally administer thrombolytic agents is necessary to accelerate clot lysis and reduce potential systemic side effects.

Medical microrobots have been proposed as the perfect candidate for such task, as they could be injected into the patient's bloodstream, navigated to the desired location and made to locally deliver a specific drug.[8,9] In the last decade, several concepts for the application of microrobots in healthcare have been proposed, such as targeted drug delivery, micro-biopsy, ablation, remote sensing and scaffolding, to name a few.[10–12] Yet, some challenges remain. For instance, with the scaling down of microrobots comes a change in the forces acting on them.[13,14] Viscous forces and surface effects become more important than volumetric effects such as inertia and weight, which pose constraints on the generation and storage of power for propulsion.[15–17] Conventional technologies used with larger robots are yet to be scaled down to micro sizes,[15] propelling researchers to look for novel methods beyond those used in conventional robotics. Magnetic manipulation has gained attention in this context, as it allows for precise steering of magnetic objects to specific locations, can wirelessly transmit power to objects, is externally controlled and biocompatible.[18] Several experimental studies using magnetic manipulation systems have been reported, considering different object geometries (e.g. sphere, rod, and helix),[19–22] locomotion strategies (i.e. propulsion, rolling, and tumbling),[20,22] and magnetic stimuli (i.e. field gradients, rotating fields, and oscillating fields).[23–25] Thus, magnetically manipulated microrobots could potentially be navigated through the cerebral arteries using external magnetic field gradients, to deliver therapeutic agents directly to the affected regions (e.g. a clot).

The ability to manipulate microrobots injected in the bloodstream depends on different parameters, such as the microrobots size, the blood flow velocity, and the magnetic field gradients imposed. Most works in the literature include the use of magnets placed near the vessels to capture magnetic microrobots. Kenjereš and colleagues investigated the deposition of microrobots (0.25–4 µm) in a brain vascular system, and concluded that the use of a magnet in the vicinity of the arteries significantly increased the microrobots capture in the targeted regions.[26] Cherry and colleagues studied the effect of microrobot size, blood velocity, and magnetic field gradients on the trapping of microrobots (15–50 µm) in a tube exposed to magnetic field gradients. They observed more trapping with larger microrobot sizes, lower flow velocities, and higher magnetic field gradients.[27] Bose and Banerjee simulated the 2D trajectory and capture of microrobots (0.25–4 µm) in a stenosed artery bifurcation to assess the influence of microrobot size, Reynolds number of the flow, and magnetic field gradients on the capture efficiency.[28] They observed that higher microrobot sizes, lower Reynolds number and higher magnetic field gradients led to higher capture efficiency. Larimi and co-workers studied the trajectory of nanomicrorobots (5–500 nm) in a bifurcation with a magnet close to one of its walls.[29] They found that the use of magnets helped to deliver microrobots to the target, though delivery decreased with increasing flow Reynolds number. Finally, Pourmehran and colleagues studied the delivery of magnetic microrobots through the human tracheobronchial airways and found that delivery to the target was higher with larger microrobot sizes and lower inhalation flow rates.[30]

Despite these contributions providing qualitative insight on the influence of various parameters over the magnetic manipulation of microrobots, they are still insufficient to enable the prediction of the magnetic field gradients needed to navigate microrobots along specific paths. Importantly, they focus on the manipulation of very small microrobots (<50 µm) and provide no quantitative information on how the magnetic gradients should be changed to navigate larger microrobots (50<Ø<1000 µm), whose higher drug-loading capacity would make them more relevant for drug-delivery applications.[31–33] Larger microrobots require lower magnetic gradients,[34] which can be induced by magnets placed farther apart. This allows for the widening of the workspaces with controlled magnetic gradients, something key for enabling the manipulation of microrobots in workspaces compatible with the human neurovascular network. Yet, the magnetic gradients needed to manipulate large microrobots in such workspaces cannot be easily inferred from the results of the above contributions.

Here we study the motion of microrobots in the blood flowing through representative cerebral bifurcations, to investigate how magnetic gradients can be used to navigate the microrobots from the injection point to a target location. First, we use computational fluid dynamics to numerically simulate the blood flow and predict its effect over the microrobots trajectory. Then, we examine the effect of microrobot diameter, artery geometry, blood velocity, injection point and target locations, on the magnetic gradients required to navigate the microrobots along the bifurcations. We use this data to develop maps and a predictive equation providing quantitative information on the required magnetic gradients, for each scenario. The developed maps and predictive equation offer easy-to-use information to guide the design, operation and optimization of automated magnetic navigation systems for healthcare applications.

## 2. Methods

### 2.1 Problem formulation

The motion of microrobots along bifurcations mimicking cerebral vascular structures was numerically studied to investigate how magnetic gradients can be leveraged to navigate the microrobots to one of the bifurcation branches. A microrobot injected into the bloodstream is affected by various forces (e.g. drag, buoyancy, and gravitational), the net of which will determine its motion through the vascular structures. To define the specific path to be followed by the microrobots when flowing through a bifurcation, we considered intermediate targets, positioned upstream and downstream the point of flow splitting, as well as at the desired outlet of the bifurcation. The magnetic gradients required to navigate the microrobots along the defined paths depend on the microrobot size, artery geometry, blood velocity, position of the intermediate targets, and microrobot initial position. For this reason, the values of these parameters were varied in a systematic way to study their influence over the required magnetic gradients.

To quantify the success in navigating microrobots along bifurcations we calculate the navigation success as the ratio between the number of microrobots reaching the target branch, and the total number of microrobots entering the bifurcation. Navigation success thus varies between 0% (no microrobots reaching the target) and 100% (all microrobots reaching the target).

### 2.2 Modelling Assumptions and Boundary Conditions

#### 2.2.1 Artery bifurcations

The bifurcations were modelled in 3D based on the geometries of the anterior, middle and posterior cerebral arteries, where most of the clots are observed.[35–38] The diameter and length of the main artery of the bifurcations ($D_1$ and $L$ in Figure S1, Supplementary Information) were defined based on the values reported in the literature (Table *1*),[39–43] to capture the main features of the blood flow in the mentioned arteries. The vessels were considered straight and smooth (Figure S1) to reduce the computational cost of the simulations. The diameter of the branching arteries ($D_{2\text{-}3}$, Figure S1) were computed using Murray's law,[44] which is based on the principle that the vascular system evolved to minimize the work associated to the blood flow. From Murray's law, the diameter of the main artery can be related to that of the branch arteries as $D_1^3=D_2^3+D_3^3$ (Figure S1), and thus $D_{2\text{-}3}$ can be computed as $D_1 \cdot 2^{-1/3}$ when $D_2=D_3$.

Table 1 – Diameters and lengths of the bifurcations considered in this study.

|  | $D_1$ [mm][45,46] | $D_2 = D_3$ [mm] | $L$ [mm][40–43] |
|---|---|---|---|
| **Anterior cerebral artery (ACA)** | 2.00 | 1.59 | 20 |
| **Middle cerebral artery (MCA)** | 2.40 | 1.90 | 30 |
| **Posterior cerebral artery (PCA)** | 1.80 | 1.27 | 10 |

### 2.2.2 Blood flow

Blood was modelled as an incompressible non-Newtonian fluid with shear-thinning behaviour[47] and density of 1060 kg·m$^{-3}$.[48] The apparent viscosity of the blood ($\eta$) was described by the Carreau model,[49] given its good agreement with experimental data[50] for both low and high shear rates:

$$\eta = \eta_\infty + (\eta_0 - \eta_\infty)[1 + (\lambda \dot{\gamma})^2]^{\frac{n-1}{2}} \quad (1)$$

where $\dot{\gamma}$ is the shear rate, $\eta_0$ = 0.056 Pa·s and $\eta_\infty$ = 0.00345 Pa·s are the limits of the apparent viscosity for low and high shear rates, respectively, and $\lambda$ = 3.313 s and $n$ = 0.3568 are fitting parameters.[51]

The blood flow was assumed to be laminar, given the associated low Reynolds number (50 < Re < 250), and non-pulsatile, given that steady flow conditions can be used to accurately predict the blood flow distribution[52–54] and many hemodynamic parameters (i.e. flow profile, vorticity, and helicity; more details in Supplementary Information).[55–59] The blood velocity differs substantially depending on the vessels it flows through.[47] Maximum blood velocities of 0.25-0.45, 0.35-0.55, and 0.45-0.65 m·s$^{-1}$ have been reported for the posterior, anterior and middle cerebral arteries, respectively, based on measurements on healthy individuals.[60–63] For this reason, to investigate a wide range of blood velocities, in this work we considered five different values of maximum blood velocities between 0.25 and 0.65 m·s$^{-1}$, in steps of 0.10 m·s$^{-1}$, for each of the three arteries considered. The velocity profile of the blood entering the bifurcations was assumed to be fully developed, with the flow development assumed to occur in the vascular structures upstream the bifurcation inlet. Based on this, the velocity profile at the bifurcation inlet was calculated using:[49]

$$u_f(r) = u_{f,max}\left[1 - \left(\frac{r}{R}\right)^{\frac{n+1}{n}}\right] \quad (2)$$

where $r$ is the radial position, $R$ is the artery radius and $n$ is a parameter describing the degree of non-Newtonian behaviour. The value of $n$ was empirically determined to ensure that the velocity profile introduced at the bifurcation inlet (i.e. predicted by Equation (2)) remains unchanged as the blood flows along the main artery of the bifurcation (in line with the assumption of fully developed flow). This was found to occur for $n$ = 0.89, as expected for shear-thinning fluids, for which $n < 1$.[49] No-slip condition was assumed at the walls of the arteries, and zero-gauge pressure was assumed at the outlets of the bifurcation.

### 2.2.3 Microrobots

The microrobots were modelled as discrete spherical microrobots whose motion was predicted using a Lagrangian approach, where each microrobot is represented as a point that is tracked in space and time, to plot its trajectory.[64] In line with other works in the literature,[26,29,30,65,66] a one-way approach was considered to model the interactions between the blood flow and the microrobots, according to which the blood flow was assumed to affect the motion of the microrobots (e.g. through drag force), while their volume was assumed not to affect the flow.[64] We chose a one-way approach instead of a two-way approach (in which the flow would affect the motion of the microrobots and their volume would affect the flow), given its much lower computational cost, which was crucial to enable the simulation of thousands of different scenarios in a timely manner.

The microrobots were assumed to consist of magnetic material only (i.e. iron oxide), with a magnetisation of 5x10$^5$ A·m$^{-1}$ equal to its saturation magnetisation ($M_s$),[67,68] and a density of 5200 kg·m$^{-3}$.[69] We chose a magnetisation equal to the saturation value so that the predicted magnetic gradients correspond to the minimum gradients producing a given magnetic force ($\vec{F}_{Magnetic}$):

$$\vec{F}_{Magnetic} = V_p \cdot M_s \cdot \nabla\vec{B} \quad \Leftrightarrow \quad \nabla\vec{B} = \frac{\vec{F}_{Magnetic}}{V_p \cdot M_s} \tag{3}$$

where $V_p$ is the particle volume and $\nabla\vec{B}$ is the magnetic field gradient.[70]

The microrobots were considered to have diameters ($d_p$) of 50, 100, 250, 500 and 1000 µm, to investigate a wide range of possible sizes for which the required magnetic gradients are necessarily different. Furthermore, we considered microrobots sizes much larger than those considered in the literature related to magnetic manipulation (where microrobot sizes are typically much smaller than 50 µm),[26–30] given their higher drug-loading capacity, which is relevant for drug-delivery applications.[31–33]

The microrobots were assumed to enter the bifurcation at five positions along the diameter of the main artery (Figure S2, left), namely, at its centre (position 3), at each of the artery walls (positions 1 and 5), and halfway between the two (positions 2 and 4). This is important to study the effect of the microrobots initial position on the magnetic gradients needed to navigate them along the bifurcations. The microrobots were assumed to be transported by the blood flow, thus entering the main artery with the blood velocity (Equation 2) prevailing in the mentioned radial positions (Figure S2, left).

### 2.2.4 Intermediate targets

As the blood flowing along the main artery of a bifurcation splits into each of its daughter branches (Figure 1), a portion of the blood stream flows through the desired branch, and the rest flows through the opposite branch. For a microrobot to be navigated through the desired branch, it must be on the portion of the blood stream that is on the same side as the desired branch (left green region, Figure 1), by the time the flow splits. Somewhat upstream this position, the trajectory of the microrobot must start to be corrected to avoid major collisions with the branch walls (e.g., around the apex of the bifurcation) that could delay its progression. To meet these two requirements, we considered the existence of two intermediate targets, upstream and downstream the flow splitting. The target upstream the flow splitting serves to guide the microrobot to the desired portion of the artery (left green region, Figure 1), whereas the target downstream the flow splitting serves to guide the microrobot through the desired branch without major collisions with its walls.

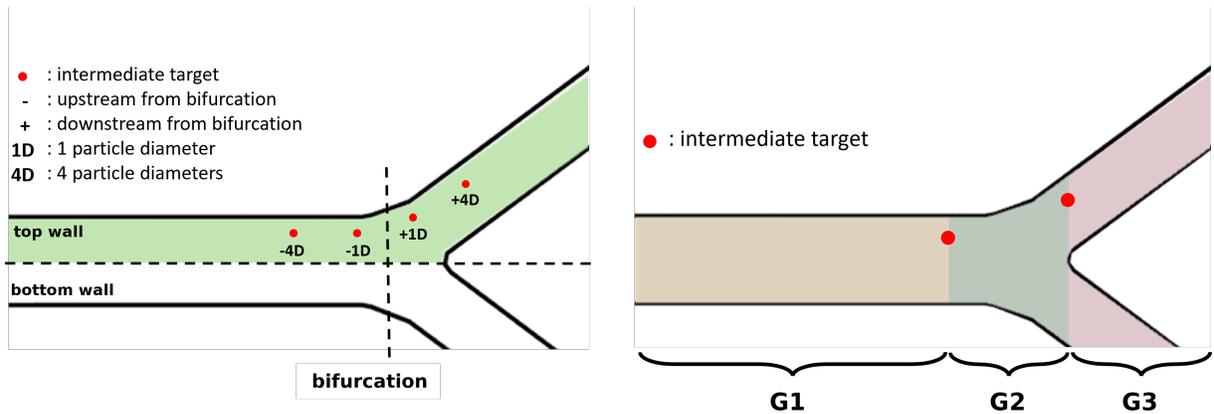

Figure 1 – Illustration of the positions of the intermediate targets considered to guide the microrobots along the bifurcation (left), highlighting the side of the branch to which the microrobots must be navigated to (green area on the left), and the segmentation of the bifurcation into three main regions defined by the position of the intermediate targets (right). Targets are positioned 1 to 4 microrobot diameters away from the position of flow splitting (bifurcation), upstream or downstream this position, to form a pair. In this illustration, 4 possible intermediate target pairings are shown: -4D+1D, -4D+4D, -1D+1D and -1D+4D (left). The flow splitting is assumed to coincide with the end of the main artery of the bifurcation. The positions of the intermediate targets define the limits of each of the three segmentations of the bifurcation, which are useful to determine the magnetic gradients G1, G2 and G3 required in each region (right).

To study how the position of the intermediate targets affects the magnetic gradients needed for navigating the microrobots along the bifurcation, we considered four different positions upstream and downstream the flow splitting, defined in relation to the microrobot diameter (i.e. -4D, -3D, -2D and -1D, for the upstream positions, and +1D, +2D, +3D and +4D, for the downstream positions, Figure 1), in a total of 16 possible combinations of target positions. Radially, the intermediate targets were positioned halfway between the centre of the arteries and the position at which the microrobots would touch the wall (Figure 1). This is important to guide the microrobots to regions with velocity lower than that at the centre (Equation (2)) – key for reducing the magnetic gradients needed for navigation – while minimizing the risk of collisions with the walls. The position of the intermediate targets will define three regions of interest (Figure 1, right) for the analysis of the required magnetic gradient, namely, the G1 region (from the inlet to the upstream target), the G2 region (between the upstream and the downstream targets) and the G3 region (from the downstream target to the outlet).

### 2.2.5 Microrobot-fluid and microrobot-wall interactions

In a classic one-way approach, the objects transported by the flow are treated as points representing the centre of the object, that can occupy any position in the domain, even impossible positions where part of the object volume falls outside the domain boundaries (e.g., artery walls, Figure S2a). To avoid this problem, in this work we have improved the calculation of the microrobot trajectories by preventing the microrobot centres from being less than one radius away from the artery walls. This constraint implies that the microrobot can touch but not cross the artery walls (Figure S2b). This is crucial to improve the accuracy of the calculated trajectories compared to the classic one-way approach and is particularly important when collisions with the walls may occur.

When computing the microrobot trajectory with this constraint, a collision occurs when the microrobot centre is one radius away from the wall, at which point the trajectory changes due to the impact. Collisions of objects with walls, and their effect over the objects velocity can be treated via the coefficient of restitution (COR), which describes the degree of elasticity of the collision. In that context, a fully elastic collision (where COR = 1) implies that the object retains all its momentum upon impact, whereas a fully inelastic collision (where COR = 0) implies the object loses all its momentum, remaining attached to the wall after impact. A survey of various numerical works involving collisions of particles with vessels walls[29,71–73] indicates that very different values of coefficients of restitution are used (0.25-1.0), without being clear which value is better for describing the situation under study in this work (e.g. the coefficient of restitution strongly depends on the mechanical properties of the objects colliding, which advises caution when trying to compare works with different underlaying assumptions). For this reason, in this work we chose a conservative value for the coefficient of restitution, i.e. COR = 1, to minimize the footprint of the underlaying assumption in the predicted trajectories, because any coefficient of restitution below unity would lead to loss of momentum and larger accumulation of microrobots near the walls, both of which would ultimately delay the progression of the microrobots along the vessels and increase the associated calculation time.

## 2.3 Numerical Methods

The blood flow in the bifurcation models was simulated using a computational fluid dynamics approach based on the finite-volume method (FVM).[74,75] Blood velocity and pressure were calculated by coupling the continuity equation, and the Navier-Stokes equation for an incompressible non-Newtonian fluid, given respectively by

$$\nabla \vec{u}_f = 0 \tag{4}$$

$$\rho_f \frac{\partial \vec{u}_f}{\partial t} = -\nabla p + \eta \nabla^2 \vec{u}_f + \rho_f \vec{g} \tag{5}$$

where $\vec{u}_f$ is the flow velocity, $\rho_f$ is the fluid density, $p$ is pressure, $\nabla$ and $\nabla^2$ are the divergence and Laplacian operators, and $\vec{g}$ is the gravitational acceleration. A steady-state, double-precision, pressure-based solver was used considering second-order discretization and a velocity-pressure coupling. Flow simulations were performed using an unstructured mesh with around $3.3 \times 10^6$ cells (maximum element size of 0.050 mm), which produced mesh-independent results with low computational cost (Figure S3). A convergence criterion of $10^{-6}$ was considered for the continuity, velocity, and viscosity, with stricter criteria producing similar results. An illustration of the unstructured mesh used is shown in Figure S3.

The three-dimensional microrobot velocity and position are obtained by integration of the force balance between the microrobot inertia and the forces that act on it over discrete time steps (Equation (6)). The total force acting on a microrobot placed in a flowing fluid can have multiple contributions such as drag, gravity, buoyancy, virtual mass, pressure gradient, Saffmann's lift and Brownian motion. The virtual mass and pressure gradient forces are only relevant if $\rho_f/\rho_p$ is close to the unity[29] (in this work $\rho_f/\rho_p \approx 0.2$), and the Saffmann's lift and Brownian forces are only relevant for sub-micron microrobots.[64] Thus, only the drag, gravity and buoyancy forces, as well as the external magnetic force, are considered in this work to compute the total (net) force $\vec{F}_T$ using the force balance equation:[76]

$$m_p \frac{d\vec{u}_p}{dt} = \underbrace{m_p \frac{\vec{u}_f - \vec{u}_p}{\tau_p}}_{\vec{F}_{Drag}} + \underbrace{m_p \frac{\vec{g}(\rho_p - \rho_f)}{\rho_p}}_{\vec{F}_{Gravity} - \vec{F}_{Buoyancy}} + \underbrace{V_p \cdot \vec{\nabla B} \cdot M_s}_{\vec{F}_{Magnetic}} \qquad (6)$$

$$\underbrace{\phantom{m_p \frac{d\vec{u}_p}{dt}}}_{\vec{F}_T}$$

where $m_p$ is the microrobot mass, $\vec{u}_p$ is the microrobot velocity, $t$ is time, $\tau_p$ is the microrobot relaxation time, $\rho_p$ is the microrobot density, $V_p$ is the microrobot volume, $\vec{\nabla B}$ is the magnetic gradient and $M_s$ is the saturation magnetization. Equation (6) is used to calculate the net force acting on the microrobot, when navigated in the blood stream under the influence of an external magnetic field gradient. For each of the system coordinates, the analytical integration in respect to time ($t$) of this equation yields the microrobot velocity along its trajectory:

$$\vec{u}_p^{i+1} = \vec{u}_f^i + e^{-\frac{\Delta t}{\tau_p}}(\vec{u}_p^i - \vec{u}_f^i) - \tau_p\left(e^{-\frac{\Delta t}{\tau_p}} - 1\right)\left[\frac{\vec{g}(\rho_p - \rho_f)}{\rho_p} + \frac{V_p \cdot \vec{\nabla B} \cdot M_s}{m_p}\right] \qquad (7)$$

where $i$ is the current iteration of the calculations. As the microrobot velocity is the derivative of the microrobot position with respect to time, $\vec{u}_p = \frac{d\vec{x}_p}{dt}$, we can obtain the microrobot trajectory by replacing $\vec{u}_p$ in Equation (7) by its derivative form, followed by its analytical integration in respect to time, yielding:

$$\vec{x}_p^{i+1} = \vec{x}_p^i + \Delta t\left[\vec{u}_f^i + \tau_p\left[\frac{\vec{g}(\rho_p - \rho_f)}{\rho_p} + \frac{V_p \cdot \vec{\nabla B} \cdot M_s}{m_p}\right]\right] + \tau_p\left(1 - e^{-\frac{\Delta t}{\tau_p}}\right)\left[\vec{u}_p^i - \vec{u}_f^i - \tau_p\left[\frac{\vec{g}(\rho_p - \rho_f)}{\rho_p} + \frac{V_p \cdot \vec{\nabla B} \cdot M_s}{m_p}\right]\right] \qquad (8)$$

The microrobot trajectory is then obtained by iterating between Equation (7) and Equation (8), for a time step $\Delta t = 10^{-5} \cdot (u_p^i + u_f^i)^{-1}$. The value $10^{-5}$ was found to be adequate to accurately compute the microrobot trajectory over time, with smaller values yielding similar trajectories.

At every time step, the magnetic gradient needed to navigate the microrobot between its current position and the next target (i.e. one of the intermediate targets or the outlet) can be obtained by rearranging Equation (8) to isolate $\vec{\nabla B}$:

$$\vec{\nabla B} = \frac{\vec{d} - \vec{u}_f \cdot t - \vec{u}_p \cdot \tau_p + \vec{u}_f \cdot \tau_p + \tau_p \cdot e^{-\frac{t}{\tau_p}} \cdot (\vec{u}_p - \vec{u}_f) + \tau_p \cdot \frac{\vec{g}(\rho_p - \rho_f)}{\rho_p} \cdot \left(\tau_p - t - \tau_p \cdot e^{-\frac{t}{\tau_p}}\right)}{\frac{V_p \cdot M_s}{m_p} \cdot \tau_p \cdot \left(t - \tau_p + \tau_p \cdot e^{-\frac{t}{\tau_p}}\right)} \qquad (9)$$

where $\vec{d} = \vec{x}_p^{i+1} - \vec{x}_p^i$ is the distance to the next target, and $t = \|\vec{d}\| \cdot \|\vec{u}_f^i\|^{-1}$ is the estimated time over which the microrobot will move between its current position and the next target. In this calculation, the microrobot is assumed to move with the velocity of the blood prevailing at its position, at every time instant.

At every iteration of the microrobot trajectory calculation, the magnetic gradient $\vec{\nabla B}$ obtained by Equation (9) is used in Equation (7) to compute the microrobot velocity, and in Equation (8) to compute the next microrobot position. This iterative process continues until the microrobot reaches one of the bifurcation outlets.

## 2.4 Validation

To validate the flow prediction, the velocity profile and the fluid viscosity obtained with the present numerical approach at three different locations along the artery (i.e., inlet, halfway of the artery length and three-quarters of the artery length), were compared with the analytical solutions of Equation (1) and Equation (2),[49,77] respectively (Figure S4, left and middle plots). The numerically predicted viscosity was also compared with measurement of various experimental works[78–81] (Figure S4, middle plot). The results of Figure S4 confirm that the predicted flow is fully developed and show a very good agreement between the obtained velocity and viscosity (hollow circles in Figure S4), and the corresponding analytical solutions and experimental data. This shows that the numerical approach developed in this work is valid and can be used to predict the flow in vascular structures.

To validate the motion and magnetic manipulation of the microrobots described by Equation (6), we did two different comparisons. In the first, we numerically replicated the work of Haverkort *et al*.[82] and compared the results obtained with our modelling approach with those in the mentioned work.

For the purpose, microrobots of different sizes (0.05 to 2 µm) were injected into a 90° bended tube where blood flowed at 0.10 m/s, and a current-carrying wire, placed in the vicinity of the tube, generated a magnetic field to capture the injected microrobots. Figure S4 shows the capture efficiency predicted using our model and that reported in the mentioned work, as a function of the size of the injected microrobots. The agreement between both results indicates that the present numerical approach can be used to predict the magnetic force and the associated effect over the trajectory of microrobots injected into flows.

In the second comparison, we ran a series of *in vitro* experiments considering a bifurcation like that of the present work, for 32 different experimental conditions involving changes in flow velocity and magnetic gradient magnitudes (see Supplementary Information for more details).

For each experimental condition, a polymeric sphere with ferromagnetic behaviour and diameter of 1.4 mm was injected into a 5 mm diameter bifurcation (Figure S5) with flowing water. As in the present work, a magnetic gradient was used to steer the sphere towards a desired outlet. Four different flow velocities and eight magnetic gradient magnitudes (Table S1) were combined following a full factorial design-of-experiments approach, where every possible combination generated an independent experimental condition that was repeated 20 times, for a total of 640 experiments (i.e. 4 flow velocities × 8 magnetic gradient values × 20 repetitions). To replicate the *in vitro* experiments, a 3D model of the bifurcation was prepared, with mesh independence tests and validation of the flow calculations done as described in sections 2.3 and 2.4. The sphere motion and trajectory along the bifurcation was calculated by iterating Equations 7-8. At each iteration (and thus position of the sphere), the magnetic field and the magnetic field gradient data measured in the experiments (Table S1) were used to calculate the sphere magnetization (Figure S6) and the resulting magnetic force acting on the sphere (Equation 3). For both the experiments and the corresponding simulations, we calculated the navigation success as the ratio between the number of polymeric spheres reaching the desired outlet and the total number of spheres entering the bifurcation (for both cases considering 20 spheres). Figure S7 shows the comparison between the navigation success observed in the experiments and that predicted by our model. The agreement between both results indicates that the modelling approach considered in this work can be used to predict the magnetic force and the associated effect over the trajectory of microrobots navigated through bifurcations.

## 2.5 Cases considered for the Numerical Simulations

To investigate how magnetic gradients can be leveraged to navigate microrobots along artery bifurcations, we considered variations in various parameters affecting the microrobot trajectory. We considered five microrobot diameters (50–1000 µm), to study sizes relevant for drug-delivery applications, and that require different magnetic gradients for navigation. We considered the geometries of three cerebral arteries where most clots are observed (ACA, PCA and MCA), considering five maximum blood velocities (0.25–0.65 m·s$^{-1}$), in line with the ranges reported for healthy individuals. Moreover, we considered five entrance positions for the microrobots defined along the diameter of the inlet, as well as 16 different combinations of intermediate target positions (from -1D+1D to -4D+4D). The values for each parameter (Table 2) were combined following a full factorial design-of-experiments approach, where each possible combination of the different parameters generated a simulation scenario. This resulted in a total of 6000 simulation scenarios (i.e. 5 microrobot diameters × 3 artery geometries × 5 blood flow velocities × 5 microrobot entrance positions × 16 target positions), which allowed studying the effects of each parameter over the microrobot trajectories and magnetic gradients required for successful navigation.

Table 2 – Parameters considered in this work that were combined in a full-factorial design of experiments to produce 6000 different simulation scenarios.

| microrobot diameter [µm] | blood velocity [m·s$^{-1}$] | entrance position | upstream target | downstream target | artery |
|---|---|---|---|---|---|
| 50 | 0.25 | top wall | -1D | +1D | ACA |
| 100 | 0.35 | mid top-centre | -2D | +2D | MCA |
| 250 | 0.45 | centre | -3D | +3D | PCA |
| 500 | 0.55 | mid bot-centre | -4D | +4D | |
| 1000 | 0.65 | bottom wall | | | |

## 3. Results and Discussion

### 3.1 Microrobot Trajectory and Magnetic Gradients

We start by focussing on a simulation scenario describing the navigation of a 500 µm microrobot released from the centre position of the inlet of the anterior cerebral artery, for a maximum blood velocity of 0.45 m/s and intermediate targets located at -2D+2D relative to the position of flow splitting, to analyse the microrobot trajectory and the magnetic gradient required to navigate it through the bifurcation. Figure 2 shows the microrobot trajectory, and the required magnetic gradient decomposed into the direction of the azimuthal angle (i.e., direction on the XY-plane) and the gradient magnitude, to enable a close inspection of the effect of the navigation parameters (i.e. direction and magnitude of the magnetic gradient).

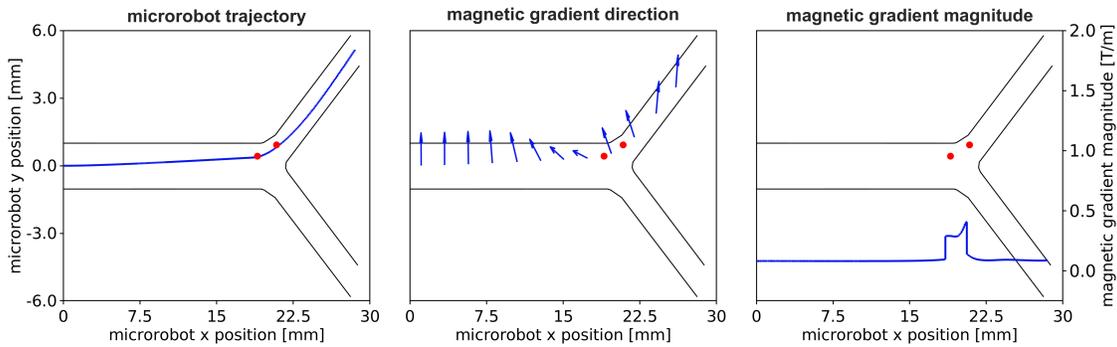

Figure 2 – Microrobot trajectory along the bifurcation (top-view of the XY-plane) and required magnetic gradient decomposed into the direction of the azimuthal angle (represented by the arrows) and gradient magnitude, for a 500 µm microrobot released from the centre position of the inlet of the anterior cerebral artery, for a maximum blood flow of 0.45 m/s and intermediate targets (red circles) placed at -2D+2D relative to the position of flow splitting.

The microrobot is successfully navigated to the desired outlet, with its trajectory following a smooth path from the inlet to the intermediate targets and to the outlet. In its motion towards the upstream target (first red circle in Figure 2), the microrobot moves along the length of the main artery (x-direction), and a small distance perpendicularly to the main flow (y-direction). While the microrobot motion along the x-direction is induced by the drag force of the blood flow, the motion along the y-direction must be induced by the application of an external force. Thus, the magnetic gradient points "upwards" in the y-direction (e.g. $x$ =5 mm in Figure 2) to induce a magnetic force that accelerates the microrobot towards the upstream target. As the microrobot moves along the bifurcation, the direction of the magnetic gradient changes to induce the deceleration (e.g. $x$ =17 mm in Figure 2), or acceleration (e.g. $x$ =24 mm in Figure 2) needed for the microrobot to reach each targets' x- and y-positions at the same time.

The plot on the right side in Figure 2 shows that the magnitude of the required magnetic gradient is very different depending on the position of the microrobots relative to the intermediate targets. The required magnetic gradients are lower between the inlet and the upstream target (G1 region, Figure 1), higher between the upstream and the downstream targets (G2 region, Figure 1), and again lower between the downstream target and the outlet (G3 region, Figure 1). This pattern of variation indicates that it is useful to analyse the required magnetic gradients separately in each of these regions of the bifurcation (i.e. regions G1, G2 and G3, Figure 1), as the associated averages are a better representation of the original data.

That the magnetic gradients are higher in the G2 region is not surprising, because that is where the microrobot trajectory must be corrected towards the second intermediate target, over a relatively small distance (i.e. the distance between the targets upstream and downstream the bifurcation). It is this dynamic variation of the magnetic gradient that allows for the desired navigation of the microrobots along the bifurcations.

### 3.2 Effect of the Parameters on the Magnetic Gradient

Here we study the influence of the different parameters on the magnetic gradients that navigate the microrobots along the artery bifurcations. The effect of each parameter on the microrobot trajectory, magnetic gradient direction and magnetic gradient magnitude is shown in Figure 3 for the anterior cerebral artery, where each parameter being considered (e.g., microrobot diameter) is varied for constant values of the others (i.e., diameter of 500 µm, blood velocity of 0.45 m/s, entrance position at the centre, and intermediate placed 2D upstream and downstream). The results of the variations are highlighted in different tints of blue, consistent across the three plots of trajectory, gradient direction, and gradient magnitude. The most relevant aspect that can be observed from Figure 3 is that the magnitude in the G1 and G3 regions remains approximately constant across all scenarios, while the magnitude in the G2 region varies significantly between parameters.

Microrobots with different diameters require different magnetic gradients to be navigated along similar trajectories (Figure 3, top charts). For instance, an increase in the microrobot diameter from 50 to 1000 µm decreases the G2 maximum magnitude from 11 to 0.5 T/m, a value more than 20 times lower, while still successfully navigating the microrobot. This decrease in the magnetic gradient magnitude is related to the effect of the microrobot diameter in the drag and magnetic forces (Equation (6)). The drag force is dependent on the microrobot surface area, whereas the magnetic force is dependent on the microrobot volume. Thus, increasing the microrobot diameter increases the relevance of magnetic forces (depending on volume) relative to the drag forces (dependent on surface area), because of the increasing volume-to-surface-area ratio of the microrobots (=$d_p/6$). This is convenient because the navigation of larger microrobots requires lower magnetic gradients and thus smaller/less powerful electromagnetic systems. Moreover, larger microrobots can carry higher doses of therapeutic agents which is interesting for drug delivery applications.[31–33] The choice of microrobot diameter should thus be weighted based on the available magnetic control system and the desired drug loading capacity.

Increasing the blood flow velocity has a very small effect over the microrobot trajectory and gradient direction, but noticeable effect over the required magnetic gradients (Figure 3, 2nd row charts). Increasing the maximum velocity from 0.25 to 0.65 m/s increased the G2 maximum magnitude from 0.15 to 0.90 T/m, a six-fold increase, because correcting the trajectory of a microrobot moving faster requires higher magnetic forces. Thus, the magnetic gradients used for microrobot navigation should be adjusted depending on the flow velocity, knowing that higher flow velocities should require stronger magnetic gradients.

The different entrance positions of the microrobot into the bifurcation only have a noticeable effect on the direction of the magnetic gradient since the trajectory of the microrobot in all five scenarios converges, as required, to the upstream target (Figure 3, 3rd row charts). As expected, the required magnetic gradient points upwards when the microrobot is released from near the bottom wall, and downwards when it is released from near the top wall. For these scenarios, the maximum G2 magnitude remains constant around 0.4 T/m. Thus, the entrance position of the microrobot does not have a meaningful impact on the magnetic gradients required to navigate the microrobots.

When varying the position of the intermediate targets (Figure 3, 4-5th row charts), as the upstream target gets closer to the point of flow division (i.e., from -4D to -1D), the maximum magnetic gradient increases slightly from 0.3 to 0.5 T/m. Similarly, as the downstream target gets farther from the point of flow division (i.e., from +1D to +4D), the maximum magnetic gradient decreases slightly, from 0.4 to 0.3 T/m. Given the different absolute positions of the different targets, the paths defined for the microrobot will have different lengths, which results in different instants at which the magnetic gradient must change the microrobot direction. Overall, the placement of the intermediate targets around the point of flow division has a small effect on the microrobot trajectory and on the magnetic gradient necessary for a successful navigation.

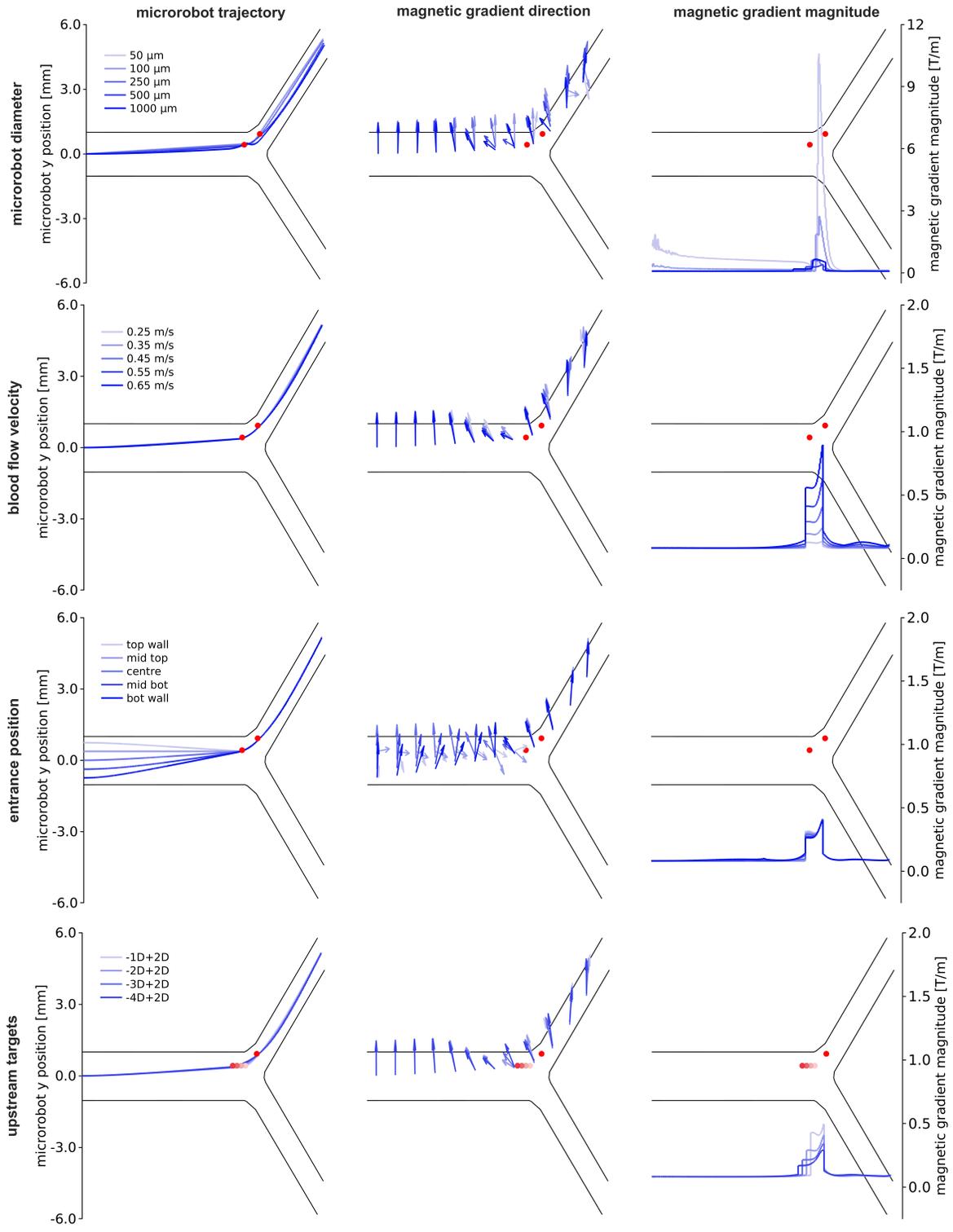

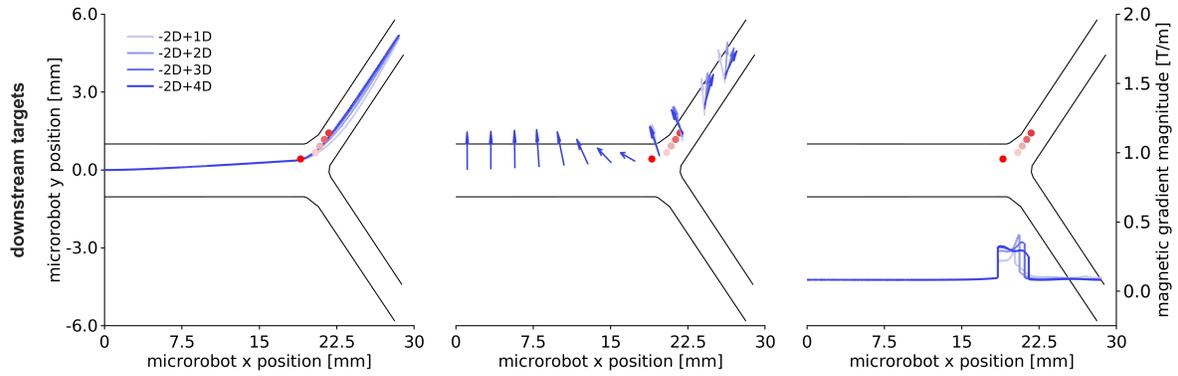

Figure 3 – Effect of microrobot diameter, blood velocity, entrance position of the microrobot, and position of the intermediate targets, on the microrobot trajectory, magnetic gradient direction, and magnetic gradient magnitude for the anterior cerebral artery model. Within each row, each parameter is changed for constant values of the other parameters (e.g., diameter of 500 µm, blood velocity of 0.45 m/s, entrance position at the centre of the artery, and intermediate targets placed 2D upstream and downstream the flow division). The red dots represent the upstream and downstream targets.

### 3.3 Maps for the Magnetic Gradients Magnitude

When computing the results of Figure 3, the required magnetic gradients were calculated considering the position and velocity of the microrobot at every time step, and the blood velocity prevailing in each position of the bifurcation. Therefore, this requires knowledge on the position and velocity of the microrobots at every time step of the calculation. Although accessing this data is straightforward in numerical simulation, that may not be the case in a real clinical scenario, where imaging techniques (e.g. angiography) would have to be used at very high frequencies, thus repeatedly exposing the patient to radiation. Yet, the fact that the magnetic gradient required to navigate the microrobots varied almost only in the G2 region, with the gradients in the G1 and G3 regions being almost constant (Figure 3), indicates that one may not need to compute the gradients at every step, because the average values in each of the mentioned G1, G2 and G3 regions may be a good representation for the "instantaneous" gradient data. For that reason, we have computed the average magnetic gradient in the G1, G2 and G3 regions, for each of the 6000 simulation scenarios considered in this work, and generated plots like those of Figure 4 (and Figures S8-S21, Supplementary Information) to show the resulting data as a function of various parameters (i.e. microrobot diameter, blood flow velocity, entrance position and target position; a spreadsheet with the raw data is given as Supplementary Information). Figure 4 shows the average gradients in the G1, G2 and G3 regions together with the number of expected collisions, as a function of the upstream and downstream target position, the microrobot entrance position and the microrobot diameter, for a given artery geometry and blood flow. This type of plot offers quantitative information on the average gradients needed to steer the microrobots along the bifurcation, and on whether wall collisions are to be expected. Based on the data in this figure, an user wishing to navigate a 1000 µm microrobot along the bifurcation with the lowest G2 gradient and number of wall collisions, would know that he should place the intermediate targets at -1D+1D or -1D+2D position (two left plots in 1$^{st}$ row of Figure 4) and use G2 gradients of 0.4–0.6 T/m. These results are in line with the results obtained in the experimental work of Belharet *et al.*,[83] which reported magnetic gradients of 0.1-0.4 T/m for a 500 µm spherical particle.

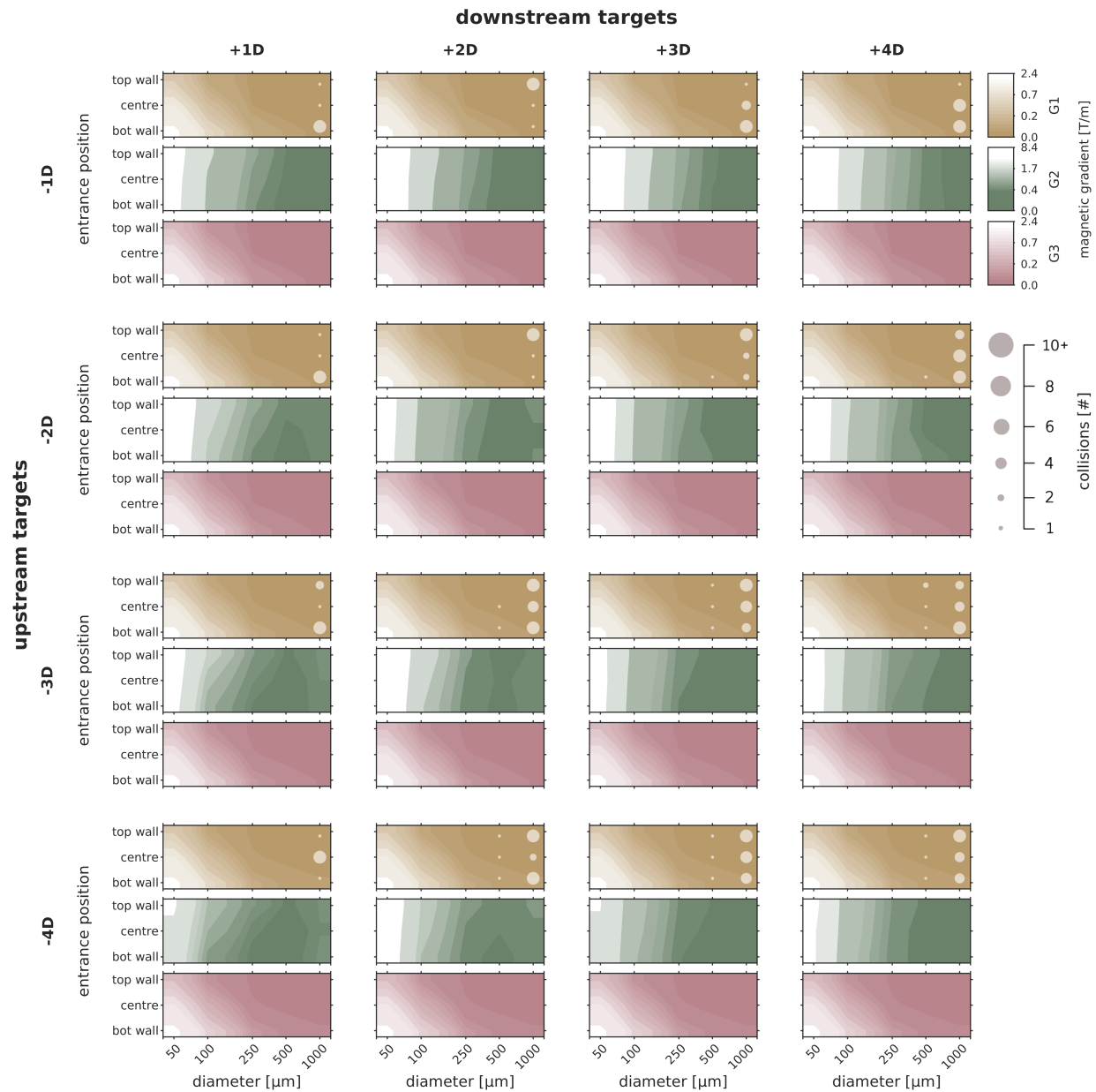

Figure 4 – Maps of the magnetic gradient magnitude in the G1, G2 and G3 regions, and number of wall collisions, for microrobots with diameters of 50, 100, 250, 500 and 1000 µm, released in the posterior cerebral artery at the top, centre and bottom positions, for intermediate target positions between -4D and +4D and maximum blood velocity of 0.45 m/s.

### 3.4 Analysis of the obtained magnetic gradients (magnitude and angle)

The maps described in the previous section compile the average magnetic gradients in the G1, G2 and G3 regions for all the 6000 simulation cases and provide useful quantitative information on the magnitude range for different scenarios. Yet, they still provide limited information on how the gradients vary across all the cases for a given parameter value (e.g., blood velocity of 0.35 m/s, or a microrobot diameter of 250 µm). To further study these variations, the average data in the G1, G2 and G3 regions were statistically analysed using boxplots, specifically regarding the magnitude and angle of the magnetic gradient. An example of this analysis is shown in Figure 5 for the G2 region, where the magnitude and angle values were grouped according to the different parameters of microrobot diameter (50–1000 µm), blood flow velocity (0.25–0.65 m/s), entrance position (top wall to bottom wall) and intermediate targets position (-1D+1D to -4D+4D). The boxplots in each row of Figure 5 were generated for each of the parameters considered by choosing a given value (e.g. microrobot diameter of 50 µm) and grouping all the data obtained for that constant parameter value (e.g. from the 6000 simulation scenarios considered, there are 1200 for which the diameter is equal to 50 µm since we analysed five different diameter values). The blue rectangle represents the interquartile range that contains 50% of the data points. The bottom and top lines of the rectangle represent the lower and upper quartiles, which represent the values below and above which 25% of the data sit. The red line inside the rectangle represents the median value, and the black horizontal lines represent the minimum / maximum values that fall 1.5x the value of the interquartile range, below / above the lower / upper quartile. Finally, the red dots represent the outliers, which are lower or higher than the minimum and maximum values, respectively.

These plots are useful to study the variability and skewness of the results of the different parameter values. For instance, when the rectangles of the boxplots are small, the interquartile range is narrow and the results have low variability, as 50% of the data points are contained between that range. Low variability suggests the gradients obtained in the different simulation scenarios are similar (meaning that the different values considered in the tested scenarios do not have a large footprint in the results). Additionally, when the interquartile range and median of two boxplots are similar, one can infer that the parameters considered in each of the two boxplots produce similar gradients.

When the diameter of the microrobot increases (Figure 5, 1$^{st}$ row plots), there is a sharp decrease in the maximum, the interquartile range, and the median value of the G2 gradient, indicating that the microrobot diameter has a large footprint in the obtained gradients. Furthermore, that the median gradient has a distinctive value for each of the tested different diameters suggests that this parameter can be used to predict the gradient required for the microrobot navigation. Additionally, that for diameters >250 µm all the results cluster around very small values, indicates low variability, thus similar magnetic gradients across all the scenarios with diameters >250 µm. On the other hand, the boxplots for the angle of the magnetic gradient in the G2 region (Figure 5, right plot in the 1$^{st}$ row) show that the results do not have distinctive median values for each of the diameters studied, with similar ranges and median values across all diameters. This shows that the G2 angle cannot be easily predicted based on the microrobot diameter, which is not surprising as the gradient direction is adjusted every time the microrobot must be accelerated / deaccelerated, as needed to change direction (between targets), or anytime it collides with the walls.

When analysing the results for the blood velocity, entrance position and target positions (Figure 5, 2$^{nd}$ – 4$^{th}$ rows), we see that the corresponding boxplots exhibit high variability, high number of outliers and no distinctive medians across the different values considered, both for the magnitude and the angle (left and right plots, Figure 5). This suggests that the median G2 magnitude and angle required for navigating the microrobot cannot be easily predicted based on the values of the blood velocity, entrance position and targets position. Additionally, the analysis of the boxplots for G2 (Figure 5) together with those for G1 and G3 (Figures S22 and S23) indicates that the median angle data appears to cluster around 90° across all tested scenarios. This was expected because an angle of 90° corresponds to the direction perpendicular to the flow along the main artery, for which the microrobot motion will depend entirely on the applied magnetic gradient (the microrobot only moves perpendicular to the flow if a magnetic gradient applies a force in that direction).

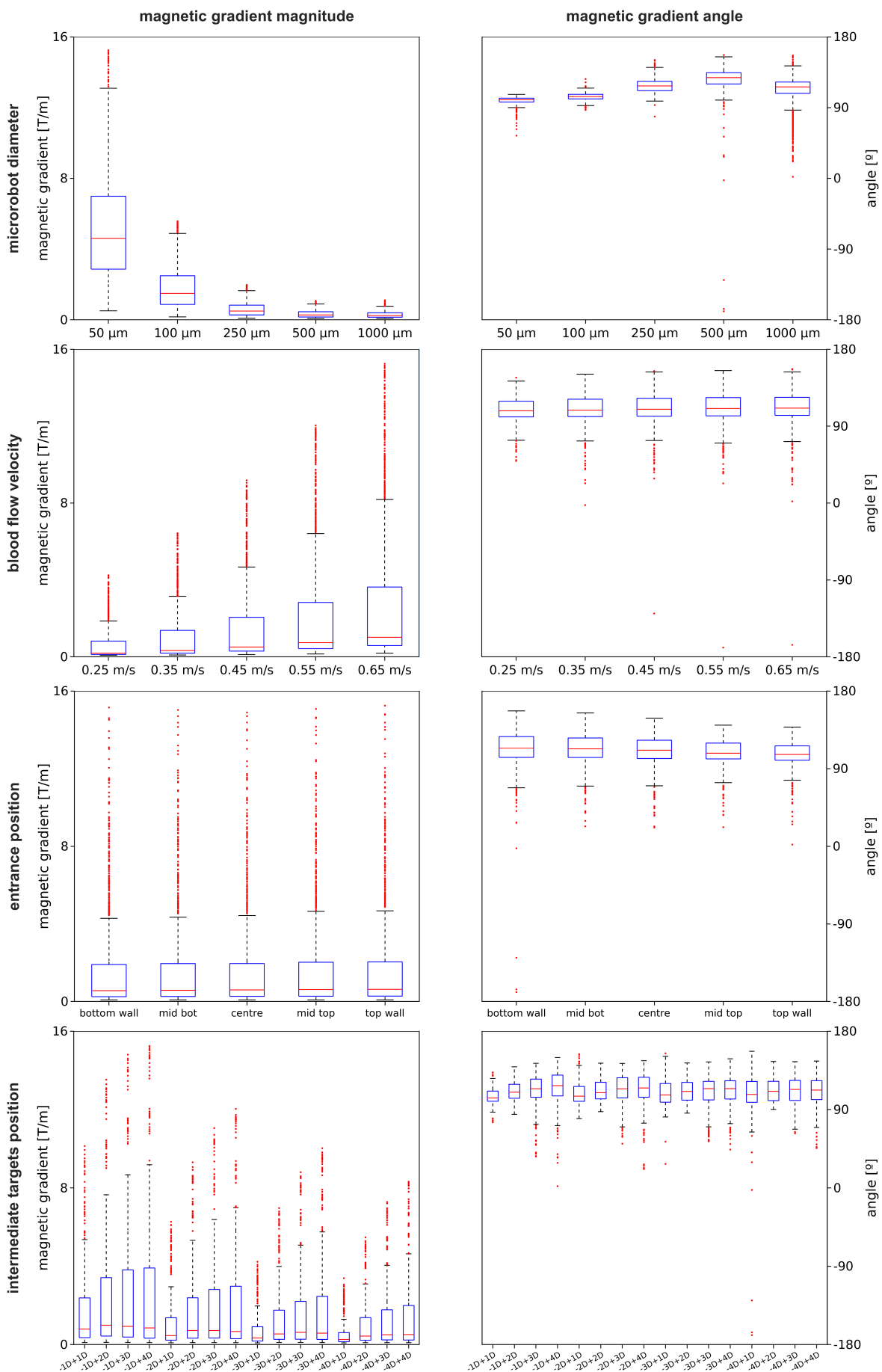

Figure 5 — Boxplots for the G2 gradient (magnitude and azimuthal angle), grouped according to the values for the microrobot diameter (50–1000 µm), blood flow velocity (0.25–0.65 m/s), intermediate targets position (-1D+1D to -4D+4D). The blue rectangle represents the interquartile range, which contains 50% of the data points; the bottom and top lines of the rectangle are the lower and upper quartiles, where 25% of the data points are below the lower quartile and 25% are above the upper quartile; the red line inside the rectangle represents the median value; the black horizontal lines represent the minimum / maximum values that fall 1.5x the value of the interquartile range, below / above the lower / upper quartile; the red dots represent the outliers, which are lower or higher than the minimum and maximum values, respectively.

## 3.5 Data-driven modelling to generate predictive equation

Having analysed in the previous section the variability of the data for the different parameters considered, we then turned our attention to the average median magnetic gradients obtained in the G1, G2 and G3 regions (Table 3). We focused on the data obtained for the minimum and maximum values of the five parameters considered in this work (i.e. microrobot diameter, blood flow, entrance position and target positions), for assessing the relative importance of changes in each of the five parameters, over the required average G1-G3 values. The data in Table 3 highlights again the importance of the microrobot diameter on the magnitude of the required gradients, with the corresponding ratio between the minimum and maximum values of G1-G3 being consistently higher than the data obtained for all the other 4 parameters (Table 3). For instance, the median G1 magnitude value for the 50 μm microrobot is 8.6 times bigger than that for the 1000 μm microrobot, and the difference for G2 and G3 is even larger, siting at 18.6 and 10.1 (Table 3, 1st row). In stark contrast, the mentioned differences for changes in blood flow, entrance position and target positions, range between 0.2 and 1.0 (Table 3, 2nd-4th rows). This shows that the microrobot diameter has the largest footprint in the obtained gradient results and can, thus, be used as an input parameter in a data-driven modelling[84] approach for generating equations that can predict the median magnitude of G1, G2 and G3, based on the chosen microrobot diameter. Such equation is useful because it should eliminate the need for running new simulations for each new set of parameters, since it could be used to predict the required gradients without running any simulation. Figure 6 shows the second order linear equations that were obtained using the least-squares method,[84] to predict the median G1, G2 and G3 values, as a function of the microrobot diameter.

Table 3 – Median values of magnetic gradients in the G1, G2 and G3 regions, obtained for the maximum and minimum values of each parameter considered (diameter, blood velocity, entrance position and target position), and ratio between the minimum and maximum values obtained.

| Parameter | minimum maximum | G1 | $\frac{G1^{minimum}}{G1^{maximum}}$ | G2 | $\frac{G2^{minimum}}{G2^{maximum}}$ | G3 | $\frac{G3^{minimum}}{G3^{maximum}}$ |
|---|---|---|---|---|---|---|---|
| Diameter [μm] | **50** <br> **1000** | 0.7069 <br> 0.0822 | 8.6 | 4.6098 <br> 0.2481 | 18.6 | 0.9064 <br> 0.0901 | 10.1 |
| Blood flow [m/s] | **0.25** <br> **0.65** | 0.0836 <br> 0.1142 | 0.7 | 0.1936 <br> 1.0169 | 0.2 | 0.0892 <br> 0.1779 | 0.5 |
| Entrance position | **Top wall** <br> **Bottom wall** | 0.0922 <br> 0.1378 | 0.7 | 0.6262 <br> 0.5603 | 1.1 | 0.1257 <br> 0.1270 | 1.0 |
| Target positions | **-4D** <br> **-1D** | 0.0933 <br> 0.0920 | 1.0 | 0.4245 <br> 0.8927 | 0.5 | 0.1269 <br> 0.1211 | 1.0 |

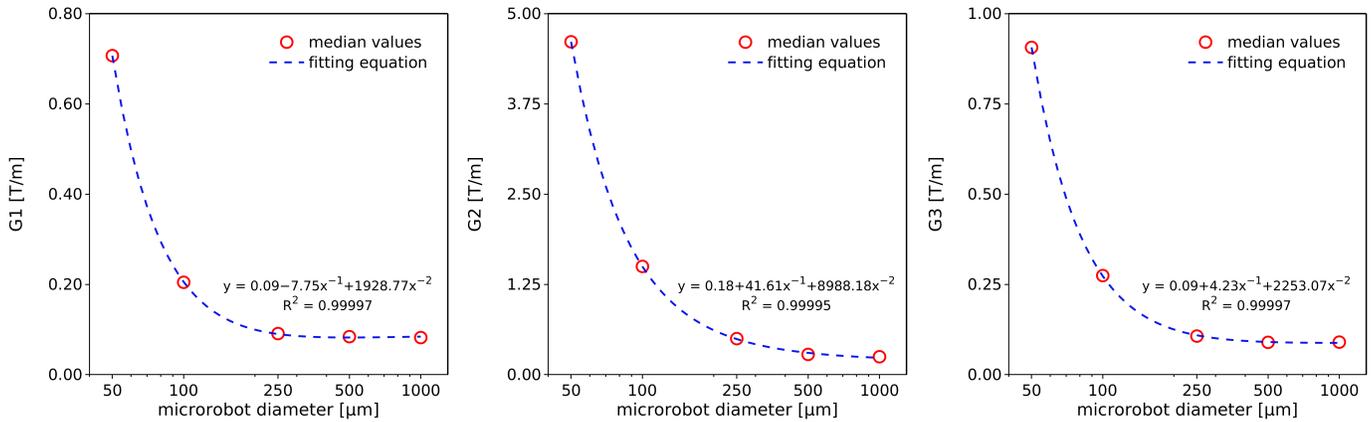

Figure 6 – Median values of G1, G2 and G3 magnitude for each microrobot diameter and respective fitting equations obtained using the least-squares method.[84] The fitting equations can be used to predict the required G1-G3 values as a function of the microrobot diameter, consisting of magnetic material only (i.e. 100% magnetic volume) with saturation magnetization. These equations can be adapted to other magnetic volume ratios and magnetic saturation values by multiplying the coefficient terms by the new magnetic volume ratio and multiplying by the ratio between the new magnetization value and the saturation value used in this work (see section 2.2.3).

After obtaining equations for predicting the required magnetic gradients in the G1, G2 and G3 regions, we checked if the predicted gradients would lead to successful navigation of the microrobot along the bifurcation. For this purpose, the entire set of 6000 simulation scenarios were re-simulated, but now considering constant median values for G1, G2 and G3, given by the developed predictive equations, instead of calculating the magnetic gradient at every time step (as done in section 3.1). Furthermore, in these new simulations, we set the gradient azimuthal angle (i.e. the angle required to navigate the microrobot towards the top outlet) and the gradient polar angle (i.e. the angle required to compensate for the effect of gravity) to a constant value of 90°, because it is the angle around which the results clustered (as explained when analysing Figure 5), but also because it is the angle at which the magnetic actuation is more efficient.

After re-running the 6000 simulations, we conducted a similar analysis to that described in section 3.4, and computed the navigation success as the ratio between the number of microrobots reaching the target branch and the total number of microrobots entering the bifurcation. When updating the magnetic gradient at every time step (dynamic method), the navigation success was 100% (meaning that, in all the 6000 scenarios, the microrobots were successfully navigated to the desired branch). When doing the navigation using the constant values of the G1-G3 gradients predicted by the developed predictive equations (constant method, Figure 6), the navigation success dropped to 94.6% (in 324 cases, the microrobot did not go to the desired branch). More importantly, the navigation was unsuccessful mostly for the smallest microrobots diameters (i.e. 50 and 100 μm), which is not particularly problematic because these diameters imply impractically large magnetic gradients, and because they are less interesting for drug delivery (because of the lower loads that can be transported in the microrobots). Furthermore, the navigation was not successful in some scenarios because the generated magnetic force was not sufficient to push the microrobot towards the top half of the bifurcation, or because collisions with the wall made the microrobot divert to the bottom outlet (Figure S24). Nevertheless, these results are very interesting because they show that the median gradients obtained for G1 G2 and G3 can be used to successfully navigate microrobots for most of the tested diameters, in a wide range of scenarios.

Having assessed the navigation success when using constant values for the magnetic gradients in the G1-G3 regions, we then focused on checking the robustness of the developed predictive/fitting equations (Figure 6). To this end, a new set of 6000 simulation cases was created using 3 new artery models (artery diameters of Table 1 increased by 0.2 mm), 5 new blood flow velocities and 5 new microrobot diameters (Table 4), and different absolute values for the entrance position and the intermediate targets (as these two parameters vary with the microrobot diameter). These 6000 new scenarios were first simulated with the dynamic method to obtain the median gradient magnitude and then simulated again using the constant median gradient magnitude calculated by the fitting equations, to assess the success of the magnetic navigation.

Table 4 - Parameter values for the new simulation scenarios

| microrobot diameter [µm] | blood flow velocity [m/s] | entrance position | upstream target | downstream target | artery |
|---|---|---|---|---|---|
| 75 | 0.30 | top wall | -1D | +1D | ACA |
| 175 | 0.40 | mid top-centre | -2D | +2D | MCA |
| 375 | 0.50 | centre | -3D | +3D | PCA |
| 650 | 0.60 | mid bot-centre | -4D | +4D | |
| 850 | 0.70 | bottom wall | | | |

We compare in Table 5 the new median G1-G3 values obtained via the dynamic method, with the G1-G3 values obtained via the fitting equations. The G1-G3 values obtained via the fitting equations underestimate the G1-G3 values obtained via the dynamic method, with average differences being 8.6% for G1, 17.1% for G2 and 10.4% for G3. The difference between the G1-G3 values obtained via the dynamic method and via the fitting equations is higher for the smaller microrobots (75 µm and 175 µm), because they imply the largest magnetic gradients due to the corresponding larger surface area to volume ratio of these microrobots (which increase the relevance of the drag force relative to the magnetic force). Nevertheless, after assessing the navigation success, we observed that the magnetic navigation was successful for 95.1% of the cases, with only 292 out of 6000 simulation scenarios having resulted in the microrobot going to the bottom outlet. This is a very positive indicator because it shows that, even with deviations between the expected and calculated values, the constant magnetic gradients obtained via the fitting equations are sufficient to navigate most microrobots to the desired branch / outlet.

Table 5 – Expected and obtained G1 / G2 values and their difference (%) for each new microrobot diameter.

| | | G1 (dynamic) | G1 (equation) | difference | G2 (dynamic) | G2 (equation) | difference | G3 (dynamic) | G3 (equation) | difference |
|---|---|---|---|---|---|---|---|---|---|---|
| microrobot diameter [μm] | 75 | 0.414 | 0.330 | 25.4% | 2.903 | 2.335 | 24.3% | 0.531 | 0.434 | 22.4% |
| | 175 | 0.118 | 0.109 | 8.2% | 0.863 | 0.714 | 20.9% | 0.159 | 0.139 | 14.2% |
| | 375 | 0.088 | 0.083 | 6.1% | 0.414 | 0.357 | 15.8% | 0.098 | 0.095 | 3.9% |
| | 650 | 0.085 | 0.083 | 2.8% | 0.297 | 0.268 | 10.9% | 0.094 | 0.089 | 5.6% |
| | 850 | 0.084 | 0.084 | 0.3% | 0.277 | 0.244 | 13.6% | 0.093 | 0.088 | 5.7% |

## 4. Conclusions

We studied the navigation of microrobots in the blood flowing through representative cerebral bifurcations, for predicting the magnetic gradients needed to navigate the microrobots until the target locations. We numerically simulated the blood flow at different velocities for three cerebral arteries, and calculated the microrobots trajectory and the magnetic gradient that would steer them along intermediate targets to the desired branch / outlet. When computing the magnetic gradients at every time step, we observed that all the microrobots were successfully navigated to the desired branch / outlet. We then analysed the magnetic gradient in three regions of the bifurcations (G1, G2 and G3), corresponding to the three paths defined by the intermediate targets. We examined the effect of microrobot diameter, artery geometry, blood velocity, injection point and target locations, on the microrobot trajectory and magnetic gradients required to navigate them along the bifurcations. We saw that an increase in the microrobot diameter from 50 to 1000 μm would decrease the maximum gradient by >20-fold. In contrast, increasing the maximum velocity from 0.25 to 0.65 m/s increased the maximum gradient by six-fold, with the entrance position of the microrobot having a noticeable effect only on the direction of the magnetic gradient. The average magnetic gradients were compiled into maps providing quantitative information on the range of values needed for different scenarios. The analysis of the results showed that the microrobot diameter was the most relevant parameter for predicting the magnetic gradients required to magnetically navigate the microrobot. We then used a data-driven modelling approach for generating second-order linear equations predicting the median magnitude of G1-G3, as a function of the microrobot diameter. We tested the robustness of the developed fitting / predictive equations and observed that microrobots could be successfully navigated along the bifurcations in the vast majority (>95%) of the tested simulation scenarios. This confirmed the robustness of the developed equations and showed the usefulness of using constant magnetic gradients dependent of the microrobot diameter, to navigate the microrobots along vascular networks. The obtained results are important for guiding the development of magnetic navigation systems, with the developed maps and predictive equations offering easy-to-use data to inform on the design, operation and optimization of such navigation systems.

## 5. Supporting Information

Rationale for considering steady blood flow conditions in the present work, description of the *in vitro* experiments of magnetic manipulation, supplementary table S1, supplementary figures S1-S24, including maps for the magnetic gradient magnitudes (S8-S21). (PDF)

Data results for generation of the magnetic gradients maps. (CSV)

## 6. Acknowledgments

Part of this work was supported by the European Union's H2020 research and innovation programme under Grant Agreements no 801464, 668786, and 952152, as well as by Horizon Europe research and innovation programme under Grant Agreement no 101047081. TSM and PGA acknowledge the support by LA/P/0045/2020 (ALiCE), UIDB/00532/2020 and UIDP/00532/2020 (CEFT), funded by Portugal through FCT/MCTES (PIDDAC). JPL also acknowledges the Agencia Estatal de Investigación (AEI) for the María de Maeztu, project no. CEX2021-001202-M, the grant PID2020-116612RB-C33 funded by MCIN/AEI/10.13039/501100011033, and the Generalitat de Catalunya for the project 2021SGR 00270.

# Supplementary Information

# Analysis of the navigation of magnetic microrobots through cerebral bifurcations

Pedro G Alves,[1,2] Maria Pinto,[3] Rosa Moreira,[3] Derick Sivakumaran,[4] Fabian Landers,[5] Maria Guix,[7] Bradley J Nelson,[5] Andreas D Flouris,[6] Salvador Pané,[5] Josep Puigmartí-Luis,[7,8,*] Tiago S Mayor[1,2,*]

[1] Transport Phenomena Research Centre (CEFT), Engineering Faculty, Porto University, Portugal
[2] Associate Laboratory in Chemical Engineering (ALICE), Engineering Faculty, Porto University, Portugal
[3] Experian, Portugal
[4] Magnebotix, Zurich, Switzerland
[5] Multi-Scale Robotics Lab, ETH Zurich, Zurich, Switzerland
[6] FAME Laboratory, Department of Exercise Science, University of Thessaly, Greece
[7] Depart. de Ciència dels Materials i Química Física and Institut de Química Teòrica i Computacional, University of Barcelona, Barcelona, Spain
[8] Institució Catalana de Recerca i Estudis Avançats (ICREA), Barcelona, Spain

* Corresponding authors: josep.puigmarti@ub.edu; tiago.sottomayor@fe.up.pt

## Steady versus pulsatile blood flow

Several studies have analysed the difference between steady and pulsatile blood flow regarding parameters such as flow characteristics (flow profile, vorticity, and helicity),[1–5] blood flow rate,[6–8] blood velocity[9] and wall-shear stress.[9,10] The major conclusion across the literature is that steady state results can approximate many hemodynamic parameters with accuracy, regarding the mean / time-average values of such parameters. For instance, Fukushima et al.[1] observed similar swirl flow in both pulsatile and steady flows in an *in vitro* symmetrical bifurcation. Malcom and Roach[5] reported that various patterns obtained in steady flow experiments on a bifurcation model were also found with pulsatile flows. Young and co-workers[11] studied the flow characteristics for a stenosed *in vitro* tube with and without flow pulsatility and concluded that many of the results obtained for various hydrodynamic factors were valid for both cases. Mikhal et al.[3] performed simulations for both steady and pulsatile flow in cerebral aneurysms and observed the same flow features in both cases, albeit with different time-dependent magnitudes. Similar findings were obtained by Karmonik and co-workers,[2] who reported similar distributions of pressure, helicity, vorticity, and velocity for steady and transient CFD simulations in cerebral aneurysms, but with different scaling factors. Chen[6] obtained similar experimental flow rates in a Circle of Willis model for steady and pulsatile flows. Hillen et al.[7] had similar findings regarding the cerebral circulation by comparing two models with and without flow pulsatility and vessel wall elasticity. The experimental results of Fahy et al.[8] showed a good prediction of in vitro blood flow distribution in the Circle of Willis using steady flow, when compared with pulsatile flow. Ku *et al.*[9] reported that the mean wall shear stress distributions observed under pulsatile flow conditions are consistent with those under steady flow, and Rabby et al.[10] observed that the wall shear stress at the beginning and end of the pulsatile flow cycles was similar, eventually becoming steady. Regarding the release and capture of particles in pulsatile flows, the work of Berselli et al.[12] shows that the particles capture region stays relatively unchanged when the particles are released at different time points over the (pulsatile) flow period. Moreover, the period-averaged capture efficiency calculated at different time instants across the flow period has very low variability. The similar capture regions combined with the low-variability of the capture efficiency show that steady, time-averaged flow conditions are good approximations to results obtained using pulsatile flows. Finally, an *in vitro* study by Bushi and co-workers[13] found that the distribution of 1.6 mm particles into two different sized daughter branches under pulsatile flow did not differ from steady flow conditions, and that the particles were distributed proportionally to the inlet flow ratios.

In summary, steady flow conditions can be used to simulate blood flow distribution in the cerebral arteries with similar results to those of pulsatile flow conditions. For this reason, in the present work, steady (non-pulsatile) blood velocities were considered in all simulations.

*In vitro* magnetic manipulation experiments

The magnetic navigation was carried out in a 8-coil electromagnetic system (OctoMag)[14] with a calibrated workspace of 6×6×6 cm in the centre. The OctoMag was used to impose a certain magnetic gradient that would steer a polymeric magnetic sphere into the desired branch of a 90° bifurcation with length of 50 mm and diameter of 5 mm (Figure S5). A constant water flow was maintained by two peristaltic pumps at different velocities (i.e. from 0.10 m/s to 0.40 m/s, in steps of 0.1 m/s). A magnetic field of 20 mT was applied in the workspace centred with the bifurcation and parallel to the direction of the flow in the main channel. Different magnetic field gradient magnitudes were imposed (i.e. from 0 to 700 mT/m, in steps of 100 mT/m), perpendicular to the direction of the flow in the main channel. The polymeric magnetic sphere had a diameter of 1.4 mm, a density of 1740 kg/m$^3$ and whose magnetization was measured for the different magnetic fields (Figure S6). Four different flow velocities and eight magnetic gradient magnitudes (Table S1) were combined following a full factorial design-of-experiments approach, where every possible combination generated an independent experimental condition. The sphere navigation along the bifurcation was repeated 20 times for each pair of flow velocity and magnetic field gradient, resulting in a total of 640 experiments (i.e. 4 flow velocities × 8 magnetic gradient values × 20 repetitions). The navigation success (Figure S7) in the experiments and in the simulations was assessed by the ratio between the number of spheres reaching the desired bifurcation outlet and the total number of spheres entering the bifurcation (for both cases considering 20 spheres).

# Supplementary Tables

Table S1 – Conditions considered for the *in vitro* experiments with magnetic manipulation of spheres navigated through bifurcations. The conditions were combined in a full factorial design-of-experiments approach, for which 20 repetitions were done, to produce 640 different experiments.

| sphere diameter [µm] | average flow velocity [m·s⁻¹] | magnetic gradient [mT/m] |
|---|---|---|
| 1400 | 0.10 | 0 |
|  | 0.20 | 100 |
|  | 0.30 | 200 |
|  | 0.40 | 300 |
|  |  | 400 |
|  |  | 500 |
|  |  | 600 |
|  |  | 700 |

# Supplementary Figures

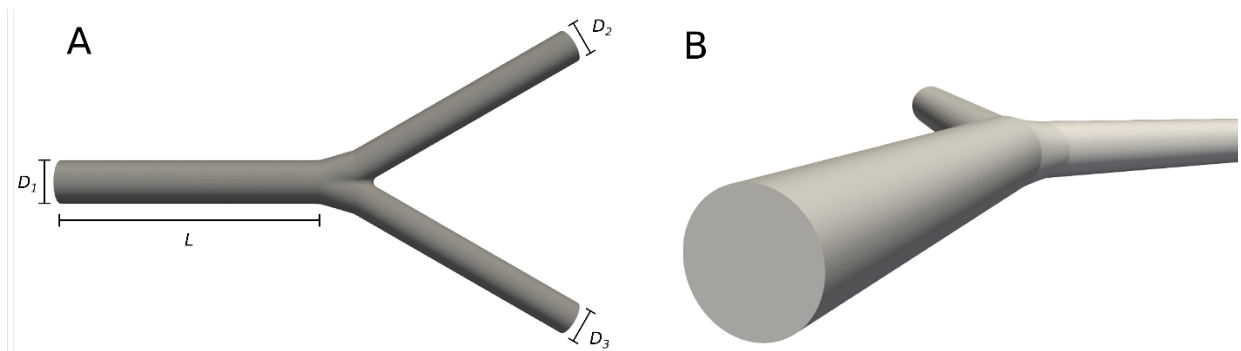

Figure S1 – Idealized 3D artery model of the anterior cerebral artery (A) viewed from the top (XY plane), and (B) viewed from the inlet. Parameters $D_1$ and $L$ are the diameter and length of the main artery, and $D_{2\text{-}3}$ are the diameters of the branch arteries.

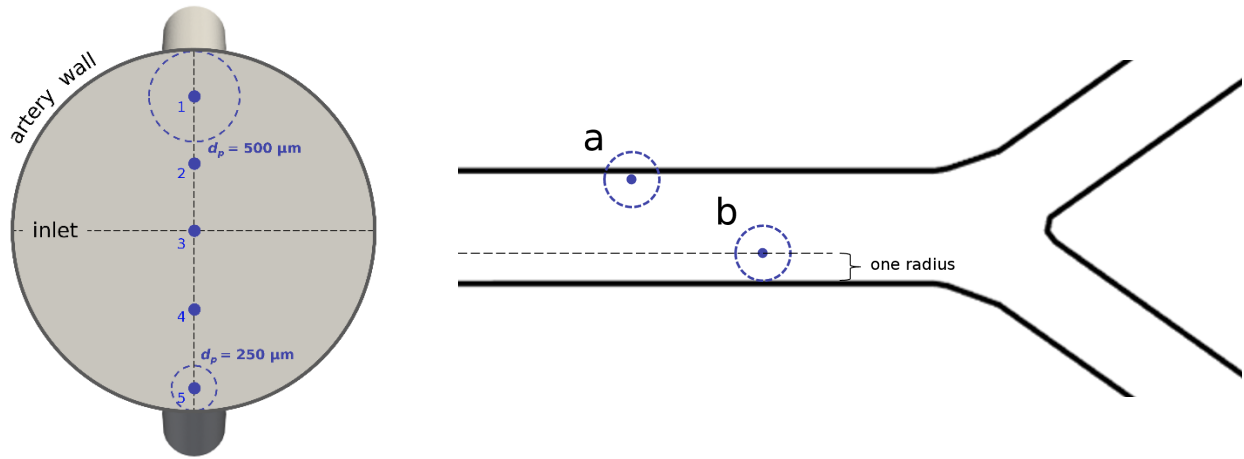

Figure S2 – Illustration of the microrobots centre at different positions in the artery, with their volume being represented with the dashed circumference. In the figure on the left are represented the entrance positions along the inlet of the main artery, for two microrobots of different diameters: 1 – top wall; 2 – halfway between top wall and centre; 3 – centre; 4 – halfway between bottom wall and centre; 5 – bottom wall. In the figure on the right, (a) represents an impossible position for a microrobot, where part of its volume falls outside the domain boundaries (i.e., the artery wall), and (b) represents the constraint imposed in this work, which prevents the microrobot centres from being less than one radius away from the artery walls to account for their volume.

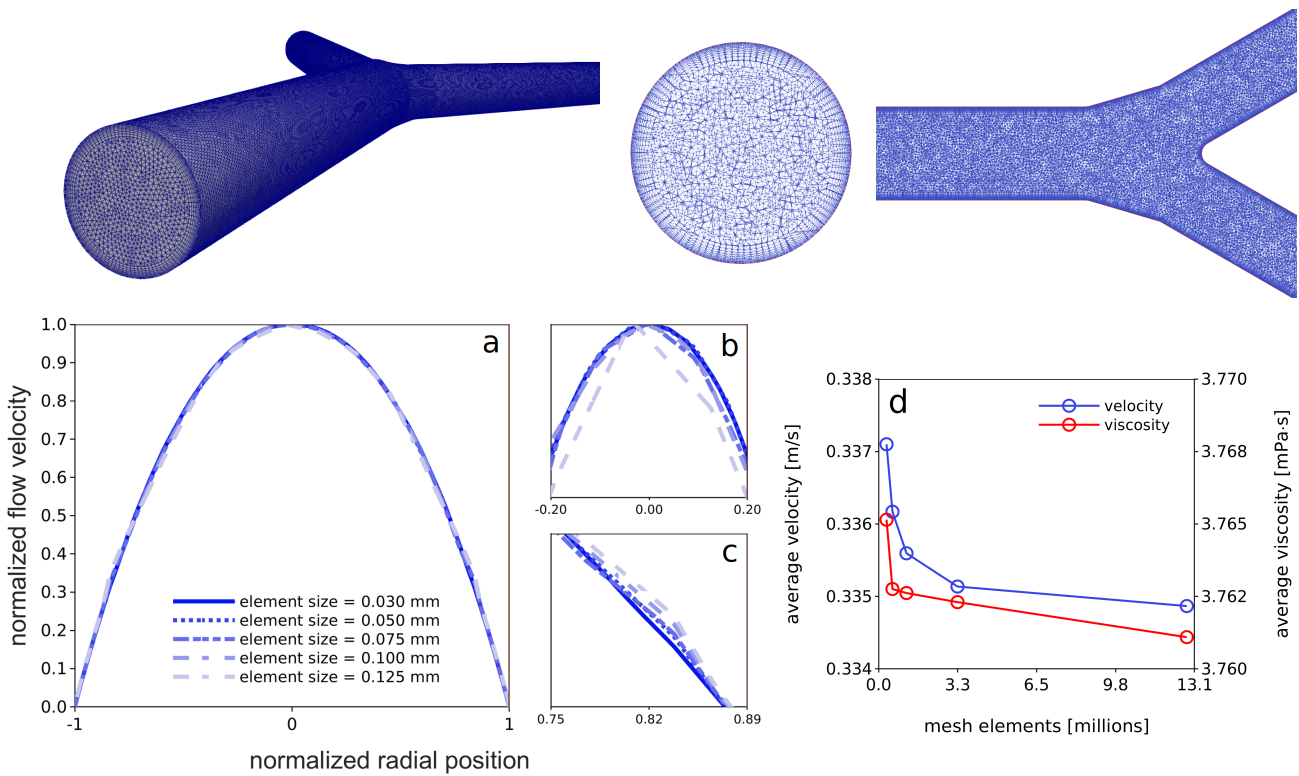

Figure S3 – Example of the unstructured mesh used in this work (complete 3D view from the inlet and mesh elements at the inlet and bifurcation area, top row), and the related mesh independence tests for meshes with five different maximum element sizes: (a) normalized flow profile along the radial position, with zoom at the (b) region of higher velocities and (c) region with higher velocity gradients; (d) average velocity and viscosity along the main artery for each mesh density (in millions of mesh elements).

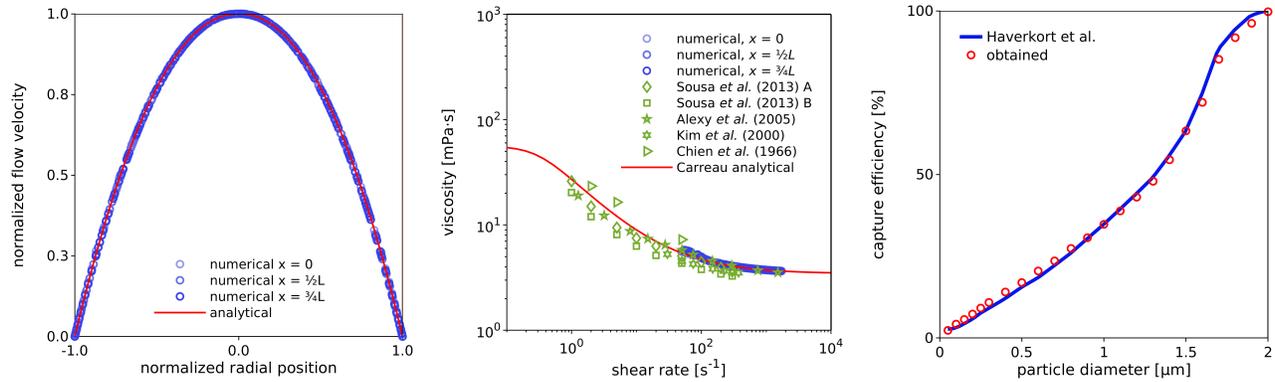

Figure S4 – Validation of the numerical approach used in this work for the flow velocity profile (left plot), blood viscosity (middle plot) and magnetic navigation of the microrobots (right plot). Comparison of predicted numerical results at different positions along the main artery (inlet, halfway of the artery length *L* and three quarters of the artery length *L*), with the analytical solutions for the flow velocity (left plot) and the blood viscosity (middle plot). The predicted numerical flow velocity is compared with the analytical solution for a fully developed flow (Equation 2), and the predicted numerical viscosity is compared with measurements of various works[15–18] and with the analytical solution obtained with the Carreau model. Comparison of the capture efficiency predicted by the present numerical approach (Equation 6) and that reported by Haverkort *et al.*[19], for microrobots of different size (right plot).

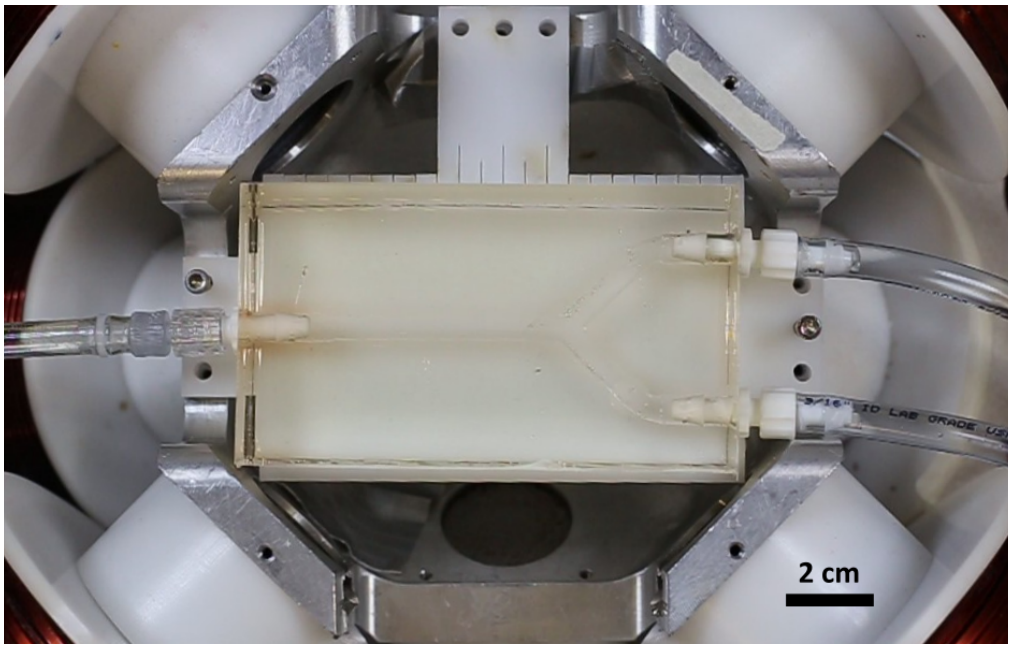

Figure S5 – Experimental setup used for the magnetic manipulation composed of an 8-coil electromagnet system (OctoMag[14]) and a 90° bifurcation with length of 50 mm and diameter of 5 mm. A constant water flow entered through the inlet on the left at a certain velocity. A single polymeric sphere was injected into the inlet and steered into the desired branch (bottom outlet on the right) by an imposed magnetic gradient. This procedure was repeated 20 times for each of the four different flow velocities (i.e. from 0.10 m/s to 0.40 m/s, in steps of 0.1 m/s) and eight magnetic field gradient magnitudes (i.e. from 0 to 700 mT/m, in steps of 100 mT/m, perpendicular to the direction of the flow in the main channel) considered, for a total of 640 experiments.

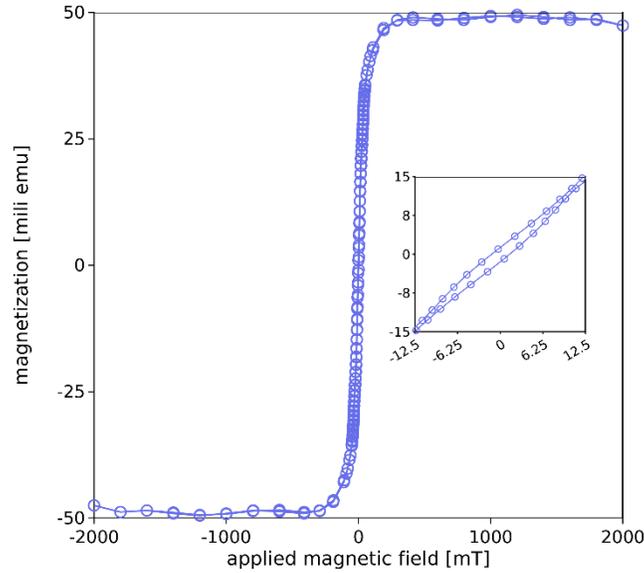

Figure S6 – Room-temperature magnetization hysteresis loop of the polymeric sphere used in the *in vitro* experiments. The inset shows the detail for the -12.5 mT to 12.5 mT range.

Figure S7 – Success of the magnetic navigation calculated the ratio between the number of spheres reaching the desired bifurcation outlet and the total number of spheres entering the bifurcation (20). The results were obtained for both the *in vitro* experiments and for the numerical simulations, considering the different conditions of imposed magnetic gradient and average flow velocity. The obtained numerical results show that the developed modelling approach allows to adequately predict the manipulation success observed in the *in vitro* experiments for the ranges of flow velocity and magnetic gradients considered. Slight differences were observed only for the lower magnetic gradients (i.e. 0 to 400 mT/m), for which the predicted number of spheres reaching the desired outlet differed by 0-4 spheres, out of the 20 considered for each experimental condition. These differences may be due to a slight mismatch between the sphere-wall restitution coefficients, and the spheres entrance position, velocity and direction, assumed in the simulations, and those prevailing in the experiments, which may have resulted in slightly different trajectories and collisions with the walls, with diversion of some of the spheres to the "wrong" bifurcation outlet.

## Maps of magnetic gradients in the G1, G2 and G3 regions

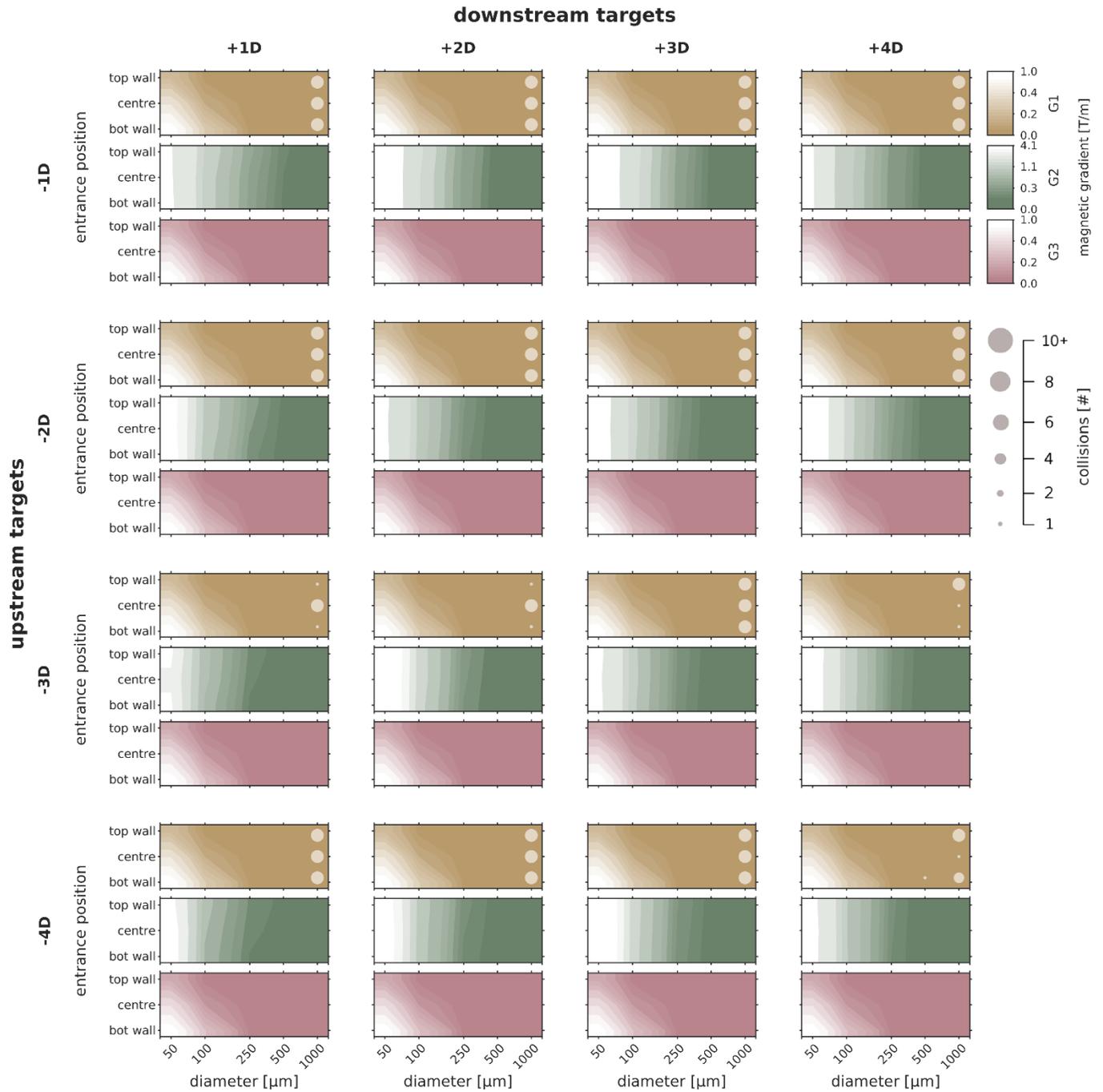

Figure S8 – Heatmaps of the magnetic gradient magnitude in the G1, G2 and G3 regions, and number of wall collisions, for microrobots with diameters of 50, 100, 250, 500 and 1000 µm, released in the anterior cerebral artery at the top, centre and bottom positions, for intermediate target positions between -4D and +4D and maximum blood velocity of 0.25 m/s.

Figure S9 – Heatmaps of the magnetic gradient magnitude in the G1, G2 and G3 regions, and number of wall collisions, for microrobots with diameters of 50, 100, 250, 500 and 1000 μm, released in the anterior cerebral artery at the top, centre and bottom positions, for intermediate target positions between -4D and +4D and maximum blood velocity of 0.35 m/s.

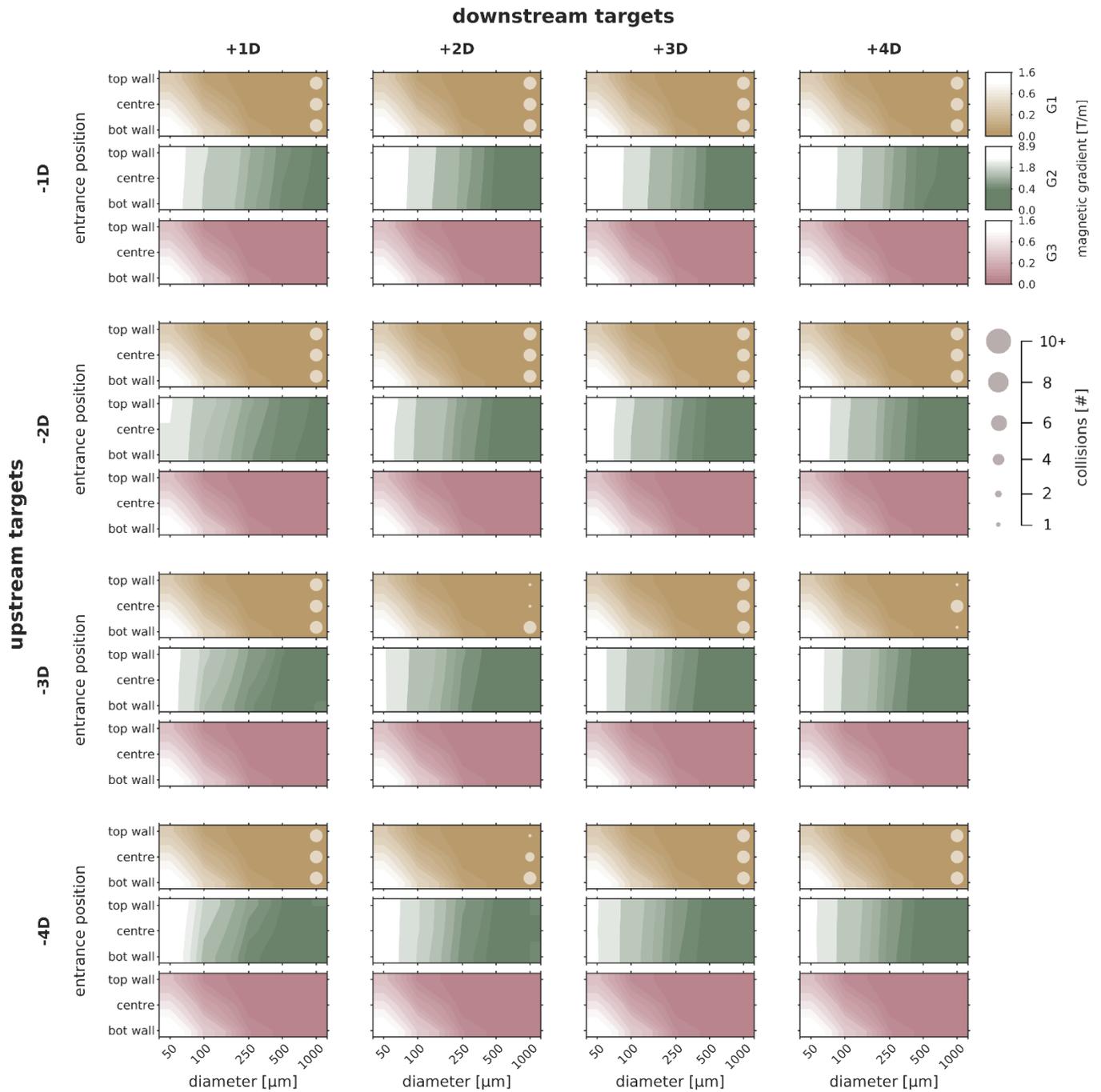

Figure S10 – Heatmaps of the magnetic gradient magnitude in the G1, G2 and G3 regions, and number of wall collisions, for microrobots with diameters of 50, 100, 250, 500 and 1000 μm, released in the anterior cerebral artery at the top, centre and bottom positions, for intermediate target positions between -4D and +4D and maximum blood velocity of 0.45 m/s.

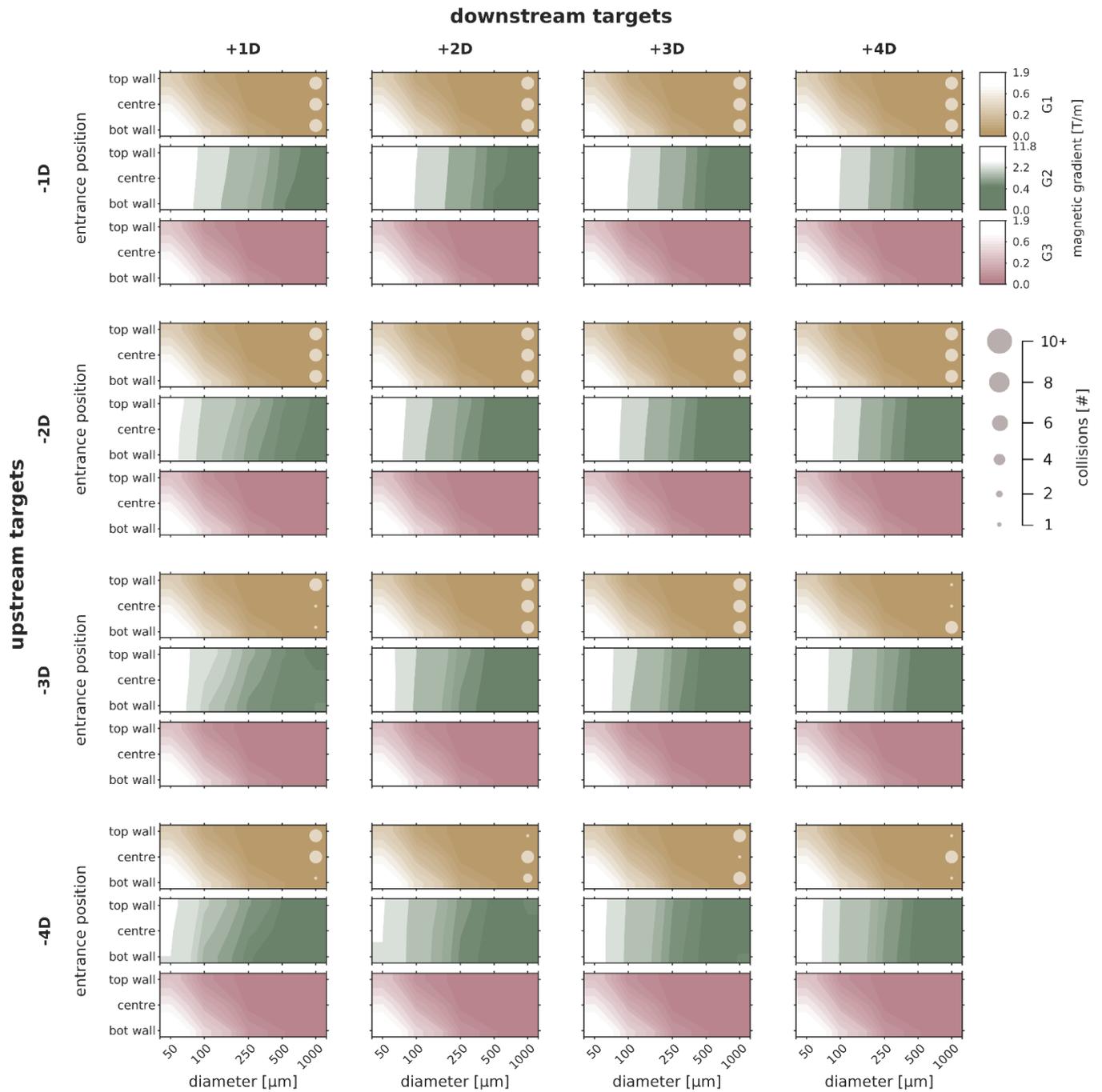

Figure S11 – Heatmaps of the magnetic gradient magnitude in the G1, G2 and G3 regions, and number of wall collisions, for microrobots with diameters of 50, 100, 250, 500 and 1000 μm, released in the anterior cerebral artery at the top, centre and bottom positions, for intermediate target positions between -4D and +4D and maximum blood velocity of 0.55 m/s.

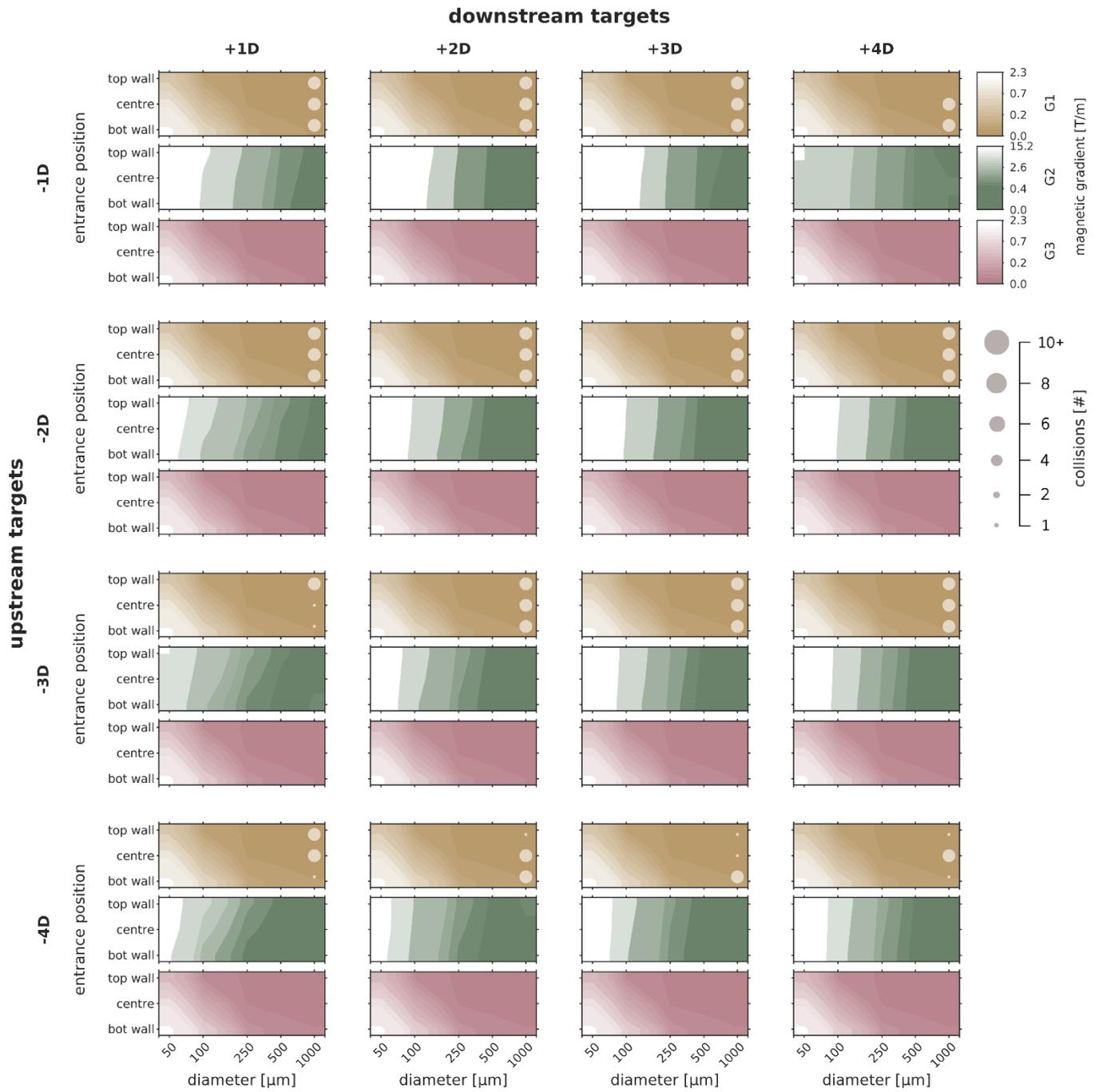

Figure S12 – Heatmaps of the magnetic gradient magnitude in the G1, G2 and G3 regions, and number of wall collisions, for microrobots with diameters of 50, 100, 250, 500 and 1000 µm, released in the anterior cerebral artery at the top, centre and bottom positions, for intermediate target positions between -4D and +4D and maximum blood velocity of 0.65 m/s.

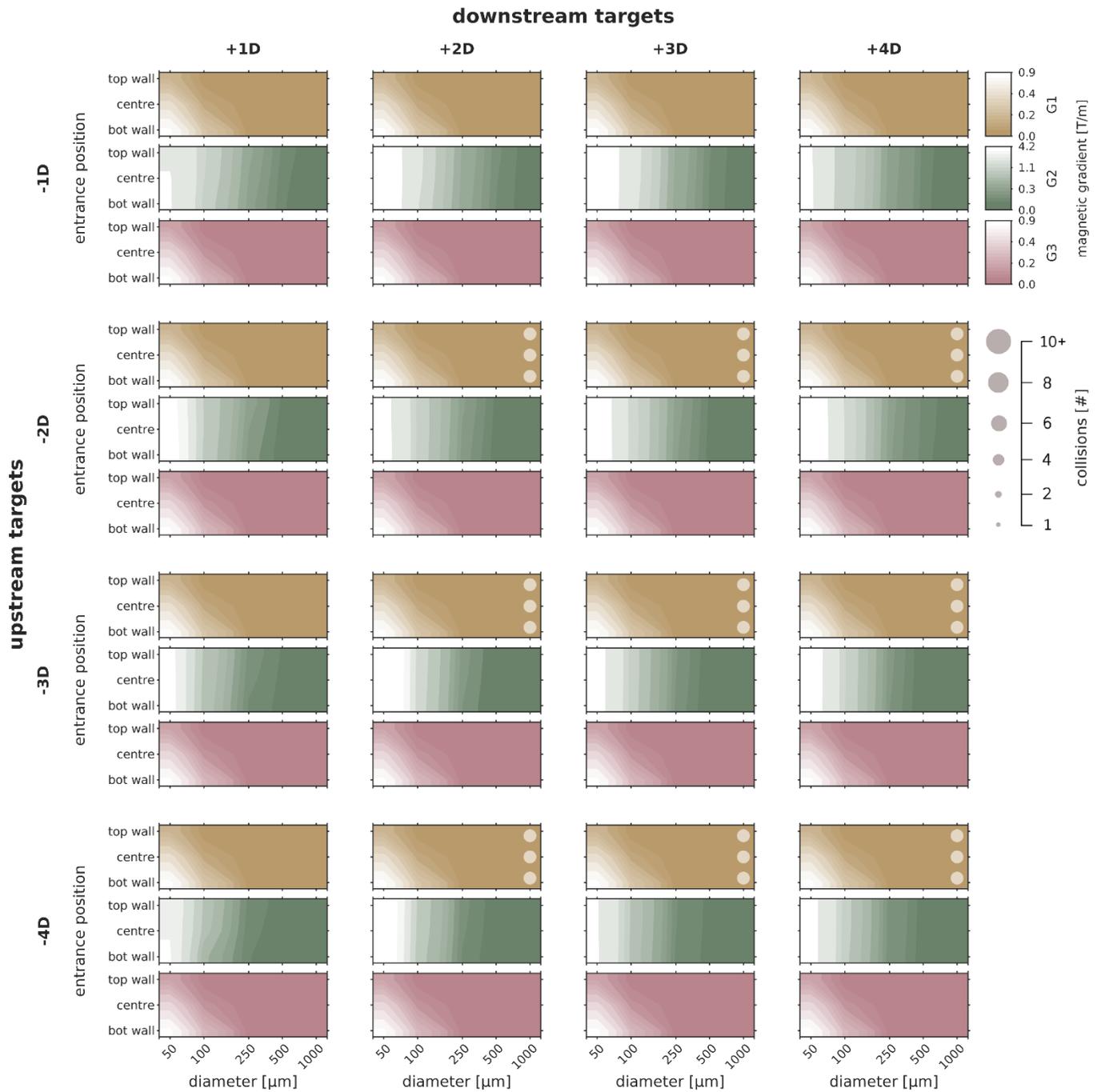

Figure S13 – Heatmaps of the magnetic gradient magnitude in the G1, G2 and G3 regions, and number of wall collisions, for microrobots with diameters of 50, 100, 250, 500 and 1000 µm, released in the middle cerebral artery at the top, centre and bottom positions, for intermediate target positions between -4D and +4D and maximum blood velocity of 0.25 m/s.

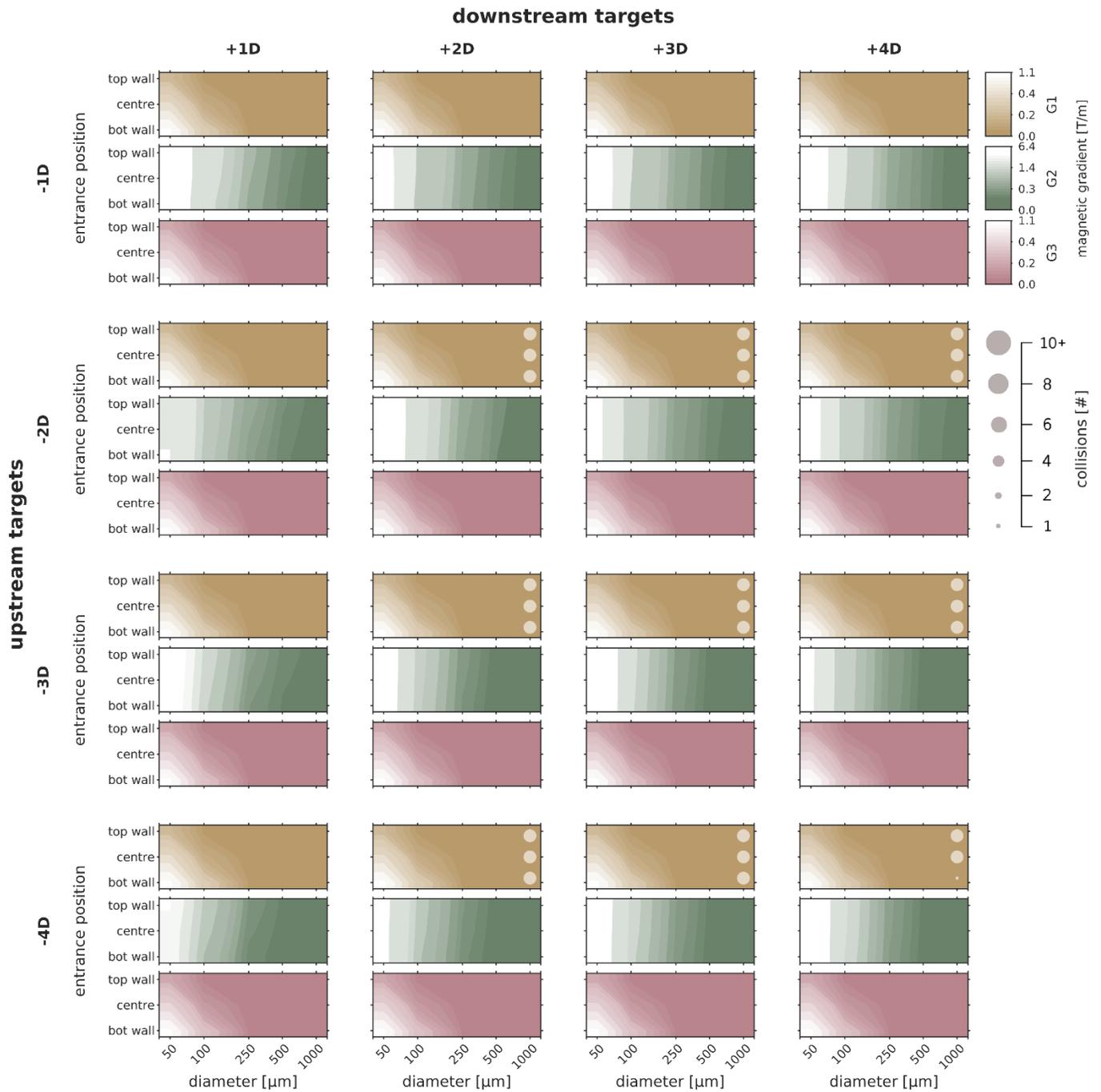

Figure S14 – Heatmaps of the magnetic gradient magnitude in the G1, G2 and G3 regions, and number of wall collisions, for microrobots with diameters of 50, 100, 250, 500 and 1000 µm, released in the middle cerebral artery at the top, centre and bottom positions, for intermediate target positions between -4D and +4D and maximum blood velocity of 0.35 m/s.

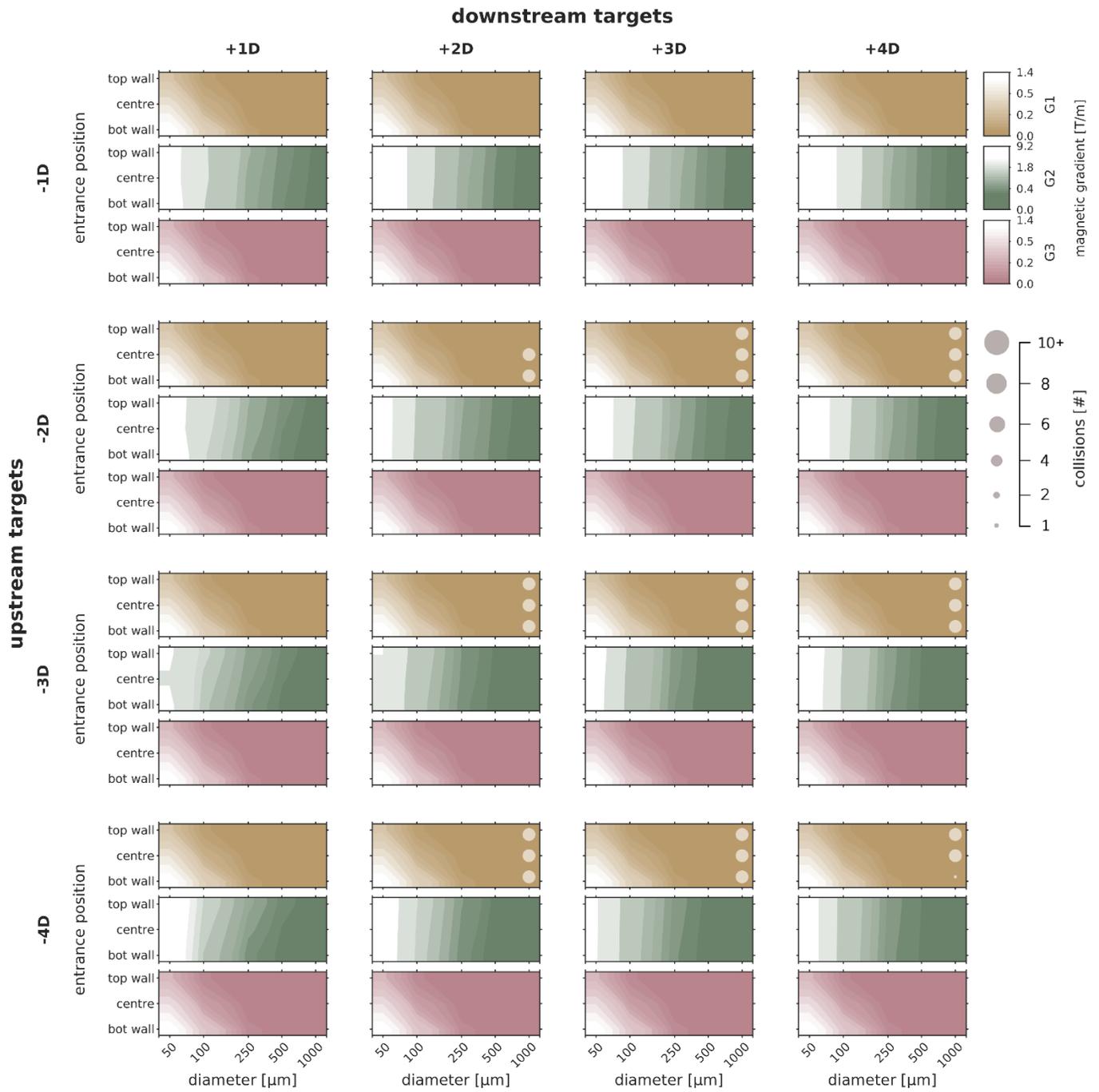

Figure S15 – Heatmaps of the magnetic gradient magnitude in the G1, G2 and G3 regions, and number of wall collisions, for microrobots with diameters of 50, 100, 250, 500 and 1000 μm, released in the middle cerebral artery at the top, centre and bottom positions, for intermediate target positions between -4D and +4D and maximum blood velocity of 0.45 m/s.

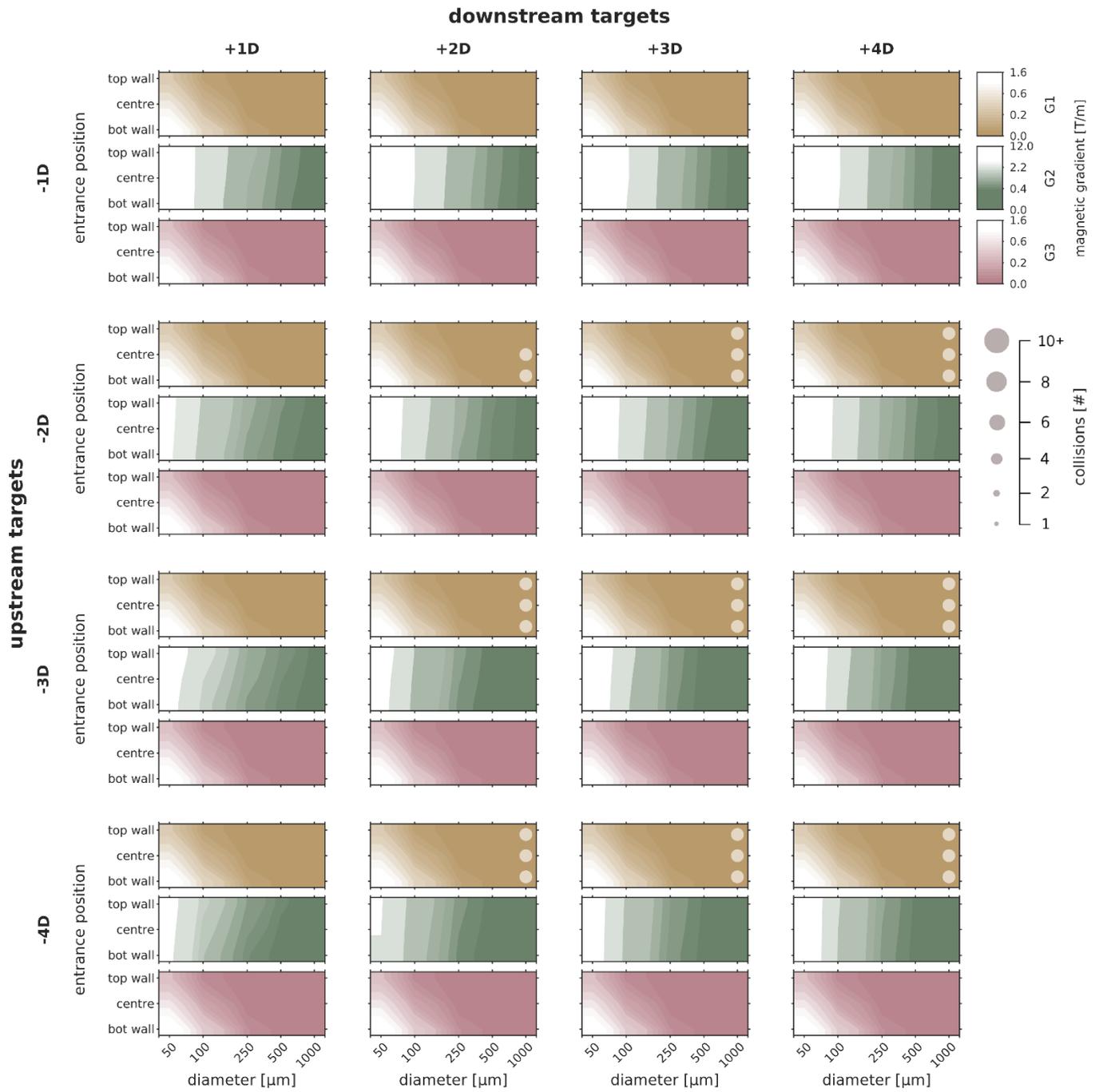

Figure S16 – Heatmaps of the magnetic gradient magnitude in the G1, G2 and G3 regions, and number of wall collisions, for microrobots with diameters of 50, 100, 250, 500 and 1000 µm, released in the middle cerebral artery at the top, centre and bottom positions, for intermediate target positions between -4D and +4D and maximum blood velocity of 0.55 m/s.

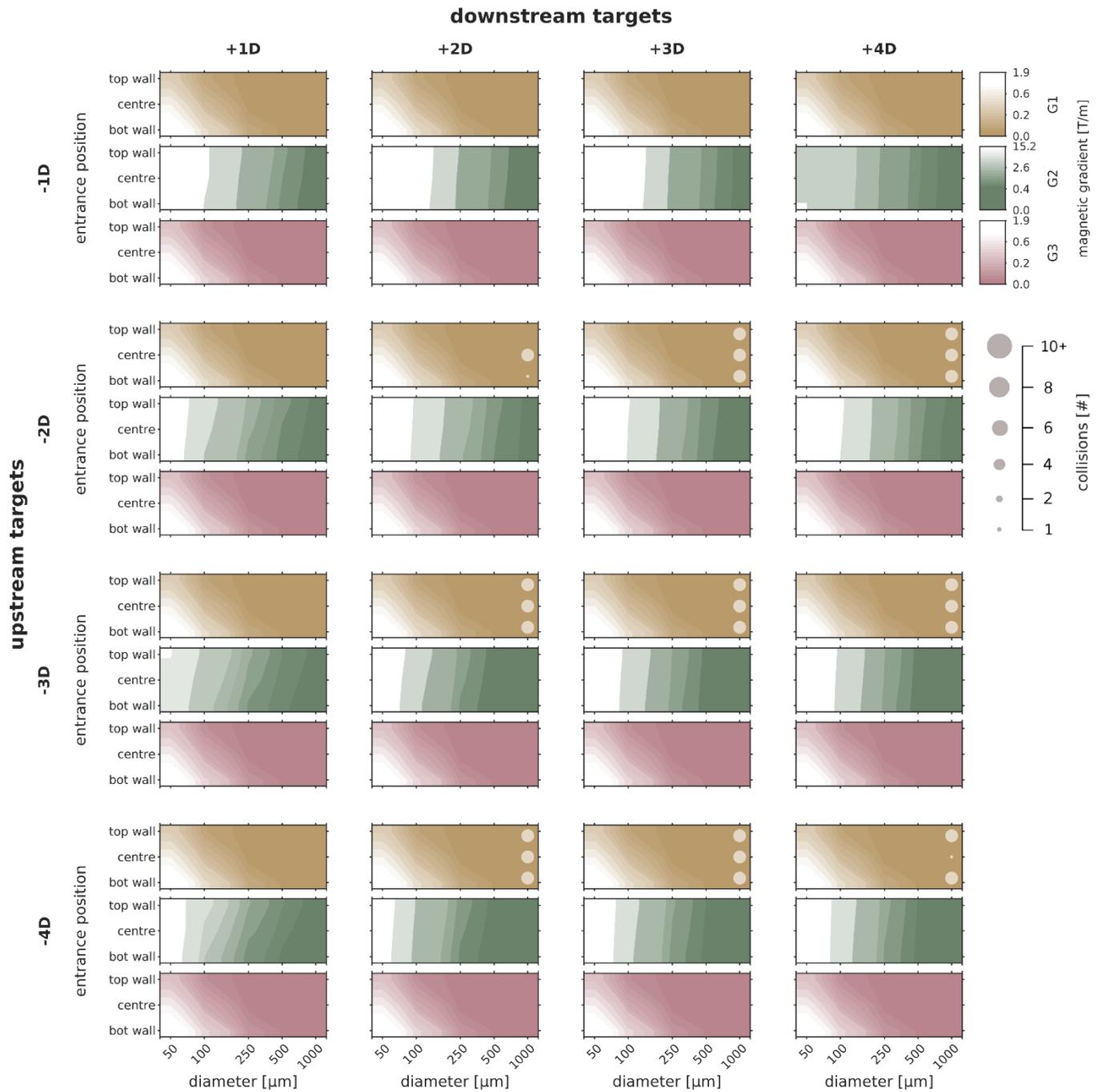

Figure S17 – Heatmaps of the magnetic gradient magnitude in the G1, G2 and G3 regions, and number of wall collisions, for microrobots with diameters of 50, 100, 250, 500 and 1000 µm, released in the middle cerebral artery at the top, centre and bottom positions, for intermediate target positions between -4D and +4D and maximum blood velocity of 0.65 m/s.

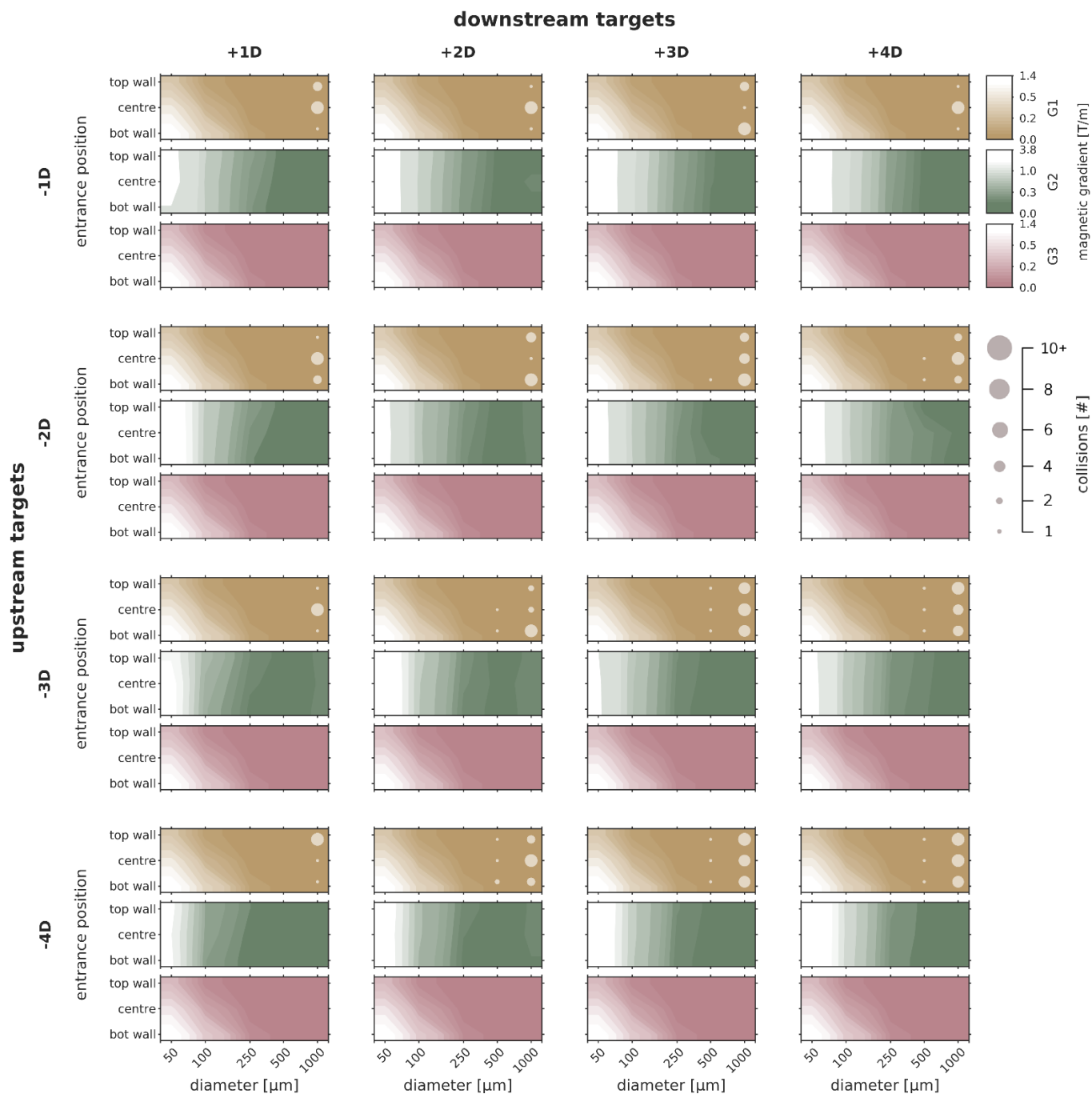

Figure S18 – Heatmaps of the magnetic gradient magnitude in the G1, G2 and G3 regions, and number of wall collisions, for microrobots with diameters of 50, 100, 250, 500 and 1000 µm, released in the posterior cerebral artery at the top, centre and bottom positions, for intermediate target positions between -4D and +4D and maximum blood velocity of 0.25 m/s.

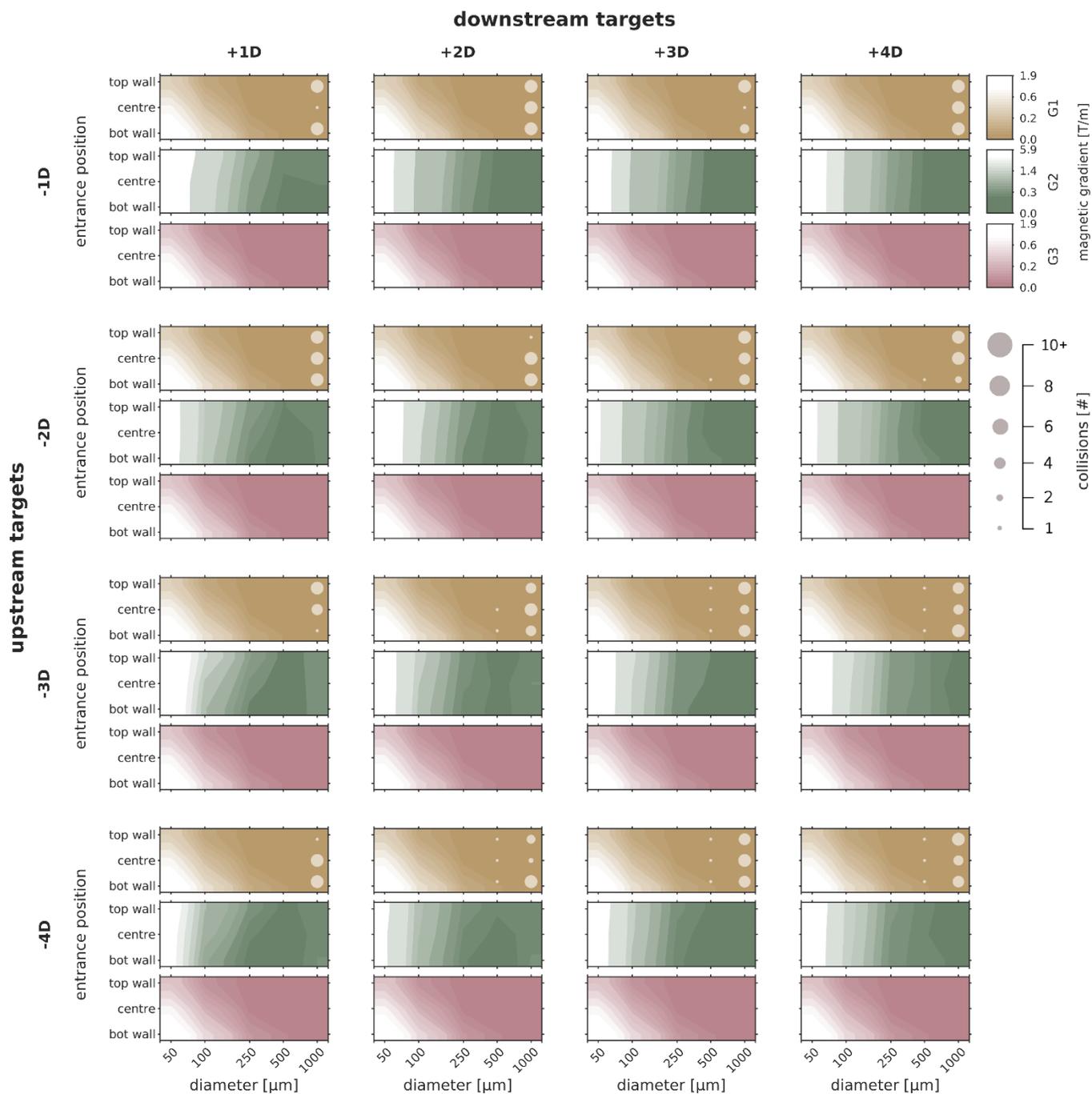

Figure S19 – Heatmaps of the magnetic gradient magnitude in the G1, G2 and G3 regions, and number of wall collisions, for microrobots with diameters of 50, 100, 250, 500 and 1000 μm, released in the posterior cerebral artery at the top, centre and bottom positions, for intermediate target positions between -4D and +4D and maximum blood velocity of 0.35 m/s.

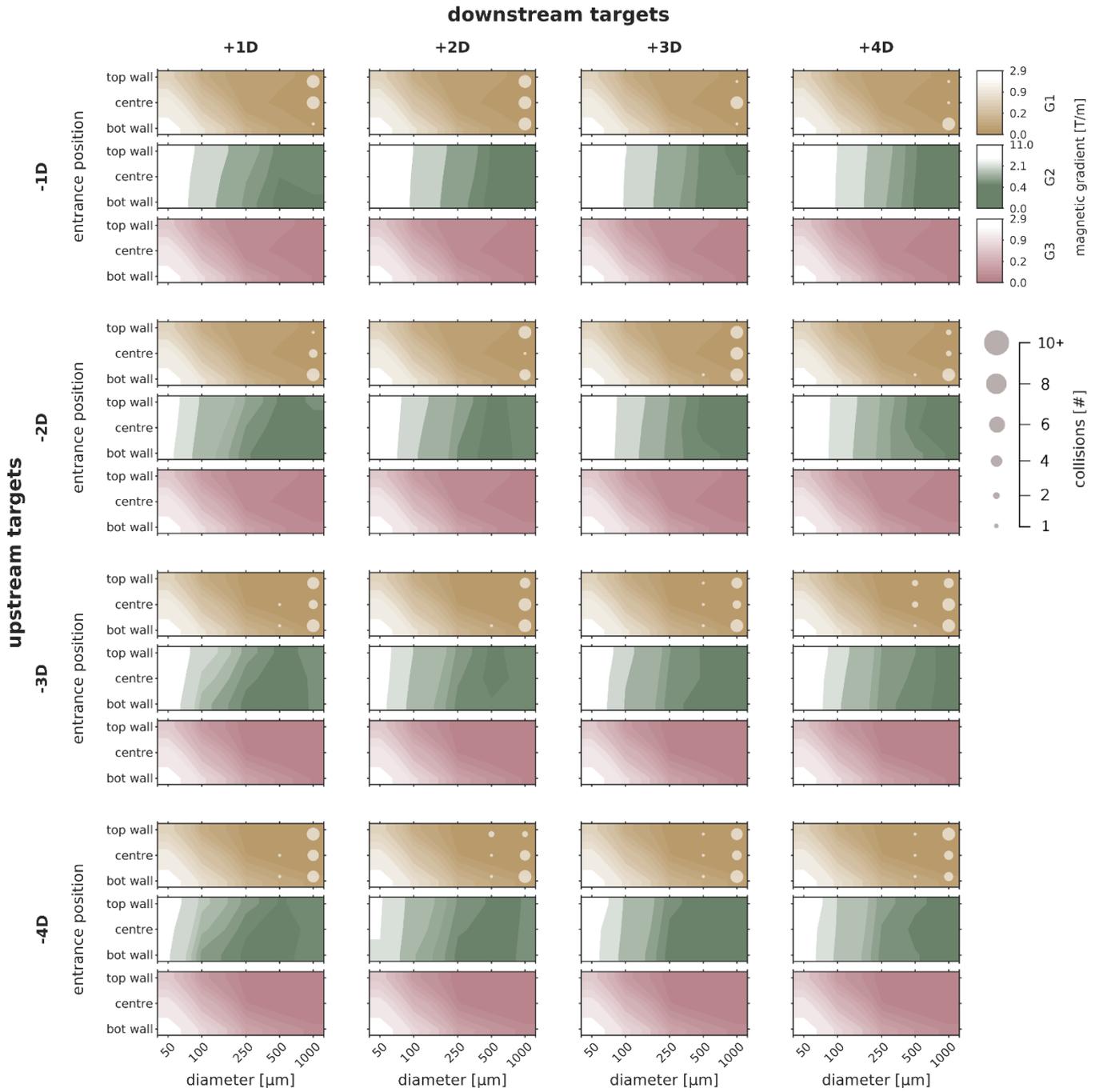

Figure S20 – Heatmaps of the magnetic gradient magnitude in the G1, G2 and G3 regions, and number of wall collisions, for microrobots with diameters of 50, 100, 250, 500 and 1000 μm, released in the posterior cerebral artery at the top, centre and bottom positions, for intermediate target positions between -4D and +4D and maximum blood velocity of 0.55 m/s.

## posterior cerebral artery  0.65 m/s

Figure S21 – Heatmaps of the magnetic gradient magnitude in the G1, G2 and G3 regions, and number of wall collisions, for microrobots with diameters of 50, 100, 250, 500 and 1000 μm, released in the posterior cerebral artery at the top, centre and bottom positions, for intermediate target positions between -4D and +4D and maximum blood velocity of 0.65 m/s.

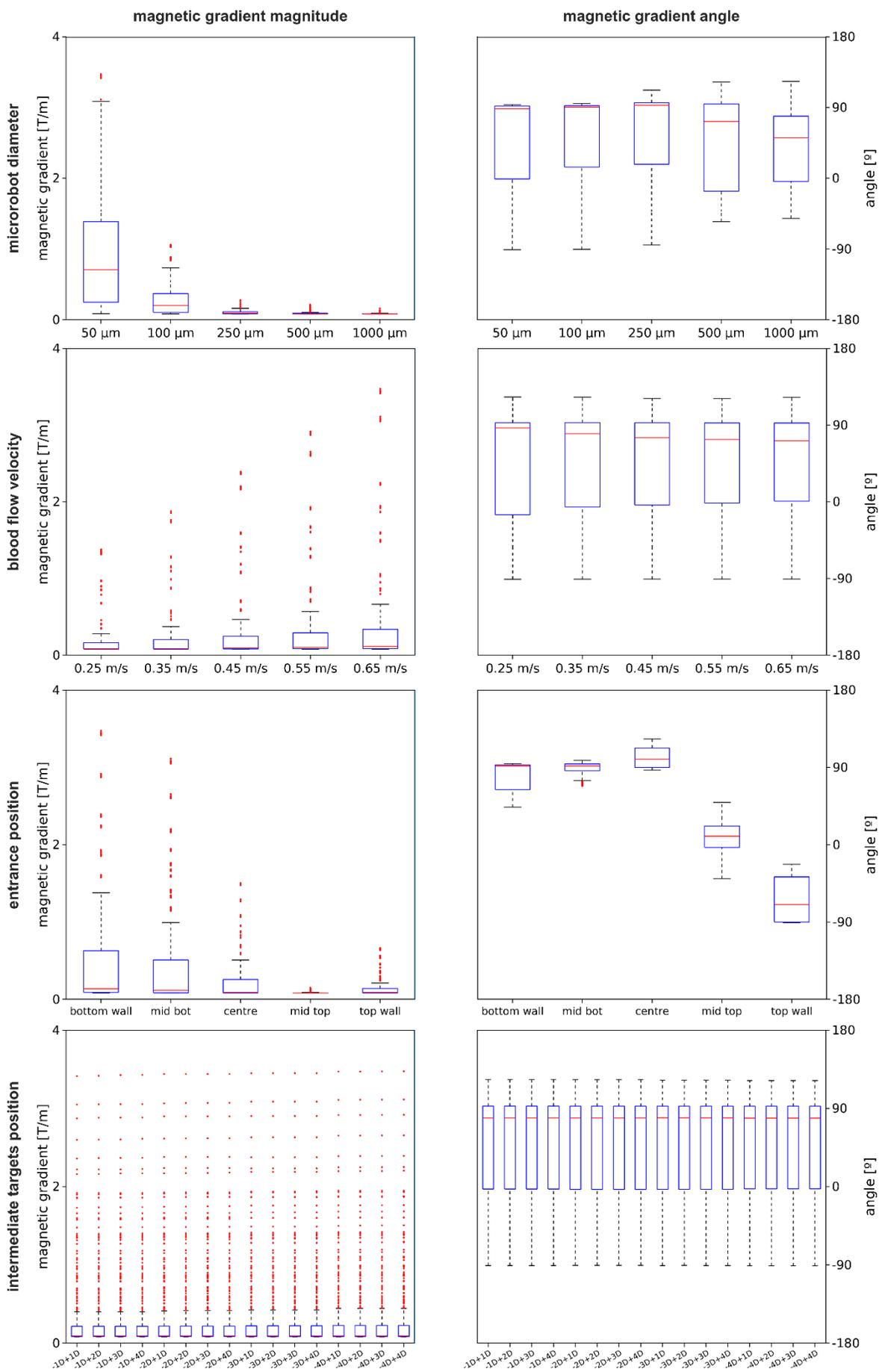

Figure S22 — Boxplots for the G1 gradient (magnitude and azimuthal angle), grouped according to the values for the microrobot diameter (50–1000 μm), blood flow velocity (0.25–0.65 m/s), intermediate targets position (-1D+1D to -4D+4D). The blue rectangle represents the interquartile range, which contains 50% of the data points; the bottom and top lines of the rectangle are the lower and upper quartiles, where 25% of the data points are below the lower quartile and 25% are above the upper quartile; the red line inside the rectangle represents the median value; the black horizontal lines represent the minimum / maximum values that fall 1.5x the value of the interquartile range, below / above the lower / upper quartile; the red dots represent the outliers, which are lower or higher than the minimum and maximum values, respectively.

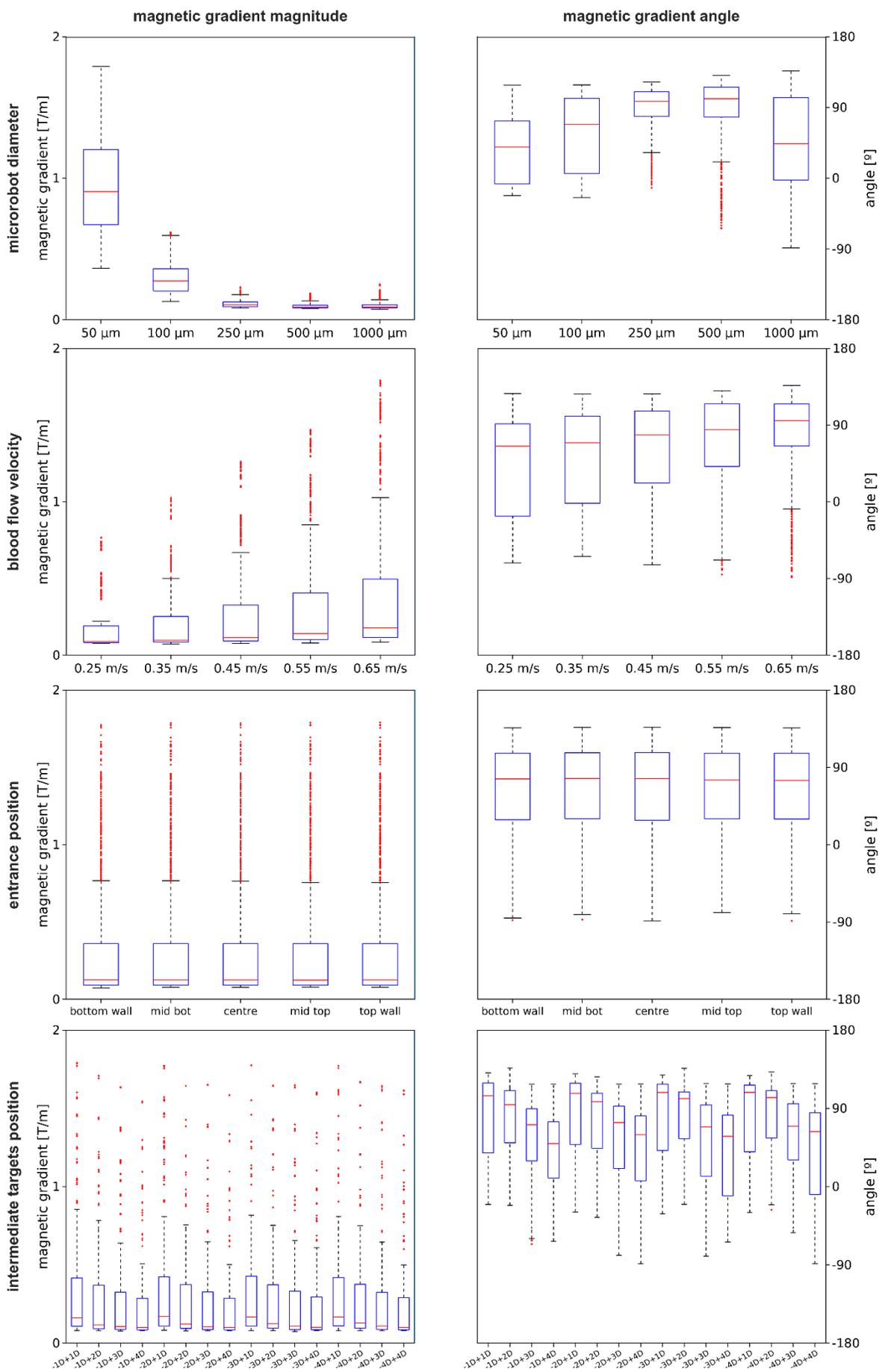

Figure S23 — Boxplots for the G3 gradient (magnitude and azimuthal angle), grouped according to the values for the microrobot diameter (50–1000 μm), blood flow velocity (0.25–0.65 m/s), intermediate targets position (-1D+1D to -4D+4D). The blue rectangle represents the interquartile range, which contains 50% of the data points; the bottom and top lines of the rectangle are the lower and upper quartiles, where 25% of the data points are below the lower quartile and 25% are above the upper quartile; the red line inside the rectangle represents the median value; the black horizontal lines represent the minimum / maximum values that fall 1.5x the value of the interquartile range, below / above the lower / upper quartile; the red dots represent the outliers, which are lower or higher than the minimum and maximum values, respectively.

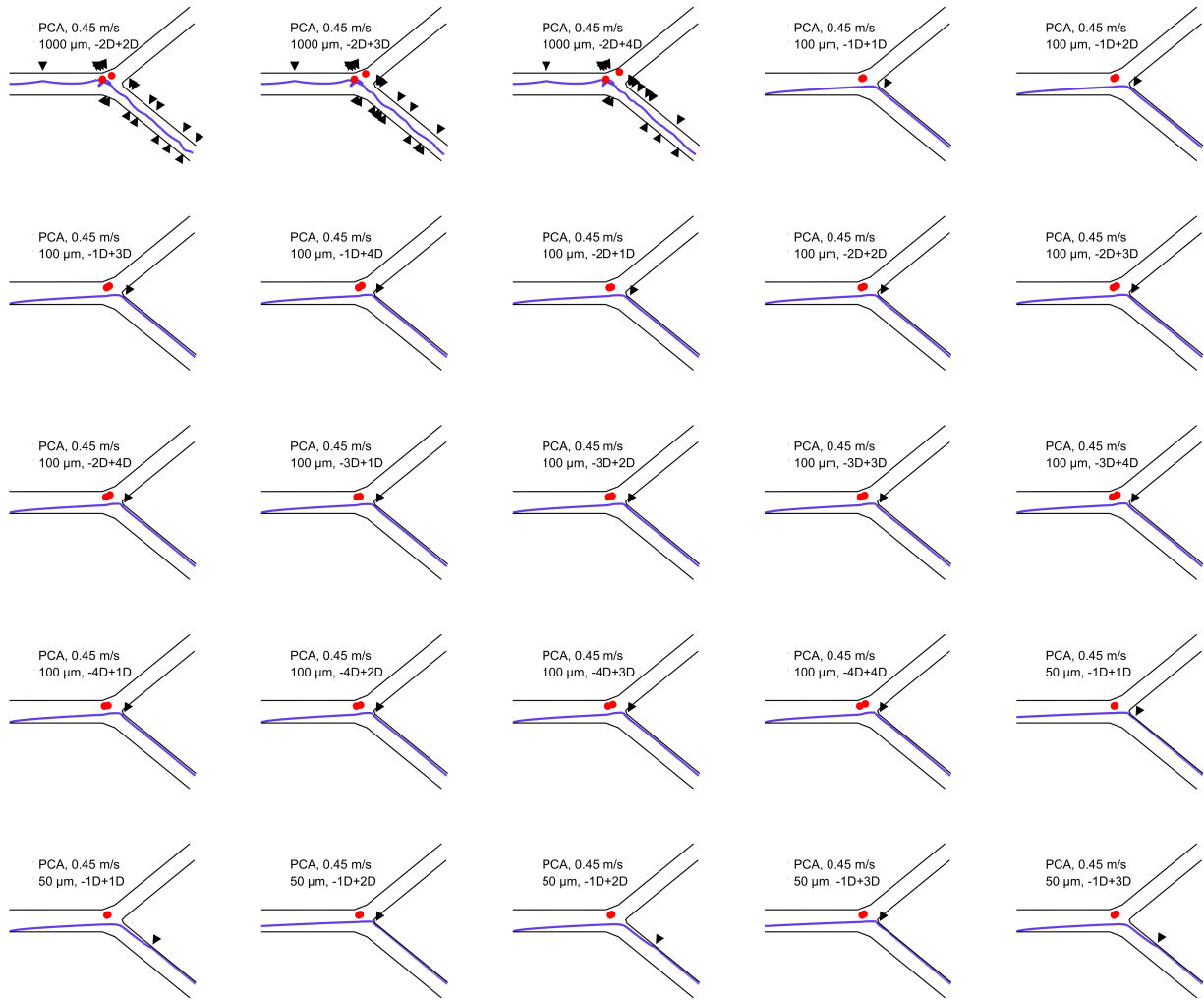

Figure S24 — Example of microrobot trajectories of unsuccessful scenarios, in which the magnetic force applied to the microrobots is not able to move the microrobot to the branch at the top half of the bifurcation, or because the collisions with the wall made the microrobot divert to the bottom outlet (three plots on the left in the first row). The intermediate targets are represented by the red dots. The position where the microrobots collide with the artery walls are marked by the black triangles.